\RequirePackage{ifpdf}
\ifpdf 
\documentclass[pdftex]{sigma}
\else
\documentclass{sigma}
\fi

\usepackage{amscd}            
\usepackage[latin1]{inputenc} 
\usepackage{eucal}            
\usepackage{bbm}              
\usepackage{stmaryrd}         
\usepackage{mathrsfs}         

\newtheorem{Theorem}{Theorem}[section]
\newtheorem{Lemma}[Theorem]{Lemma}
\newtheorem{Corollary}[Theorem]{Corollary}
\newtheorem{Proposition}[Theorem]{Proposition}
{\theoremstyle{definition}
 \newtheorem{Example}[Theorem]{Example}
\newtheorem{Definition}[Theorem]{Definition}
\newtheorem{Remark}[Theorem]{Remark}
}

\renewcommand{\mathbb}[1]{\mathbbm{#1}} 
\newcommand{\I}          {\mathrm{i}}
\newcommand{\E}          {\mathrm{e}}
\newcommand{\D}          {\operatorname{\mathrm{d}}}
\newcommand{\cc}[1]      {\overline{{#1}}}
\newcommand{\id}         {\operatorname{\mathsf{id}}}
\newcommand{\image}      {\operatorname{{\mathrm{im}}}}
\newcommand{\Lie}        {\operatorname{\mathscr{L}\!}}

\newcommand{\Pol}        {\operatorname{\mathrm{Pol}}}
\newcommand{\tr}         {\operatorname{\mathsf{tr}}}
\newcommand{\ad}         {\operatorname{\mathrm{ad}}}

\newcommand{\End}        {\operatorname{\mathsf{End}}}
\newcommand{\Der}        {\operatorname{\mathsf{Der}}}
\newcommand{\cl}         {\mathrm{cl}}
\newcommand{\divergence}  {\operatorname{\mathrm{div}}}
\newcommand{\at}[1]      {\big|_{#1}}

\newcommand{\Anti}       {\Lambda}
\newcommand{\Sym}        {\mathrm{S}}
\newcommand{\insa}       {\operatorname{\mathrm{i}_{\mathrm{a}}}}
\newcommand{\inss}       {\operatorname{\mathrm{i}_{\mathrm{s}}}}

\newcommand{\Schouten}[1]{\left\llbracket{#1}\right\rrbracket}
\newcommand{\Deg}        {\operatorname{\mathrm{Deg}}}
\newcommand{\dega}       {\operatorname{\mathrm{deg}_{\mathrm{a}}}}
\newcommand{\degs}       {\operatorname{\mathrm{deg}_{\mathrm{s}}}}
\newcommand{\degsst}     {\operatorname{{\mathrm{deg}}_\ver}}
\newcommand{\Euler}      {\operatorname{\mathsf{H}}}
\newcommand{\Eulerfib}   {\operatorname{\mathcal{H}}}
\newcommand{\ver}        {\mathsf{ver}}
\newcommand{\hor}        {\mathsf{hor}}
\newcommand{\lift}       {\mathsf{lift}}
\newcommand{\Unlift}     {\operatorname{\mathsf{P}}}
\newcommand{\varsharp}   {{\scriptscriptstyle \#}}
\newcommand{\starweyl}  {\mathbin{\star_{\scriptscriptstyle\mathrm{Weyl}}}}
\newcommand{\starstd}   {\mathbin{\star_{\scriptscriptstyle\mathrm{Std}}}}
\newcommand{\starastd}  {\mathbin{\star_{\cc{\scriptscriptstyle\mathrm{Std}}}}}
\newcommand{\starkap}     {\mathbin{\star_{\kappa}}}
\newcommand{\starStrkap}    {\mathbin{\star'_{\kappa}}}
\newcommand{\starkapStr}    {\mathbin{\star_{\kappa'}}}
\newcommand{\astkap}    {\mathbin{\ast_{\kappa}}}
\newcommand{\astStrkap}    {\mathbin{\ast'_{\kappa}}}
\newcommand{\astkapStr}    {\mathbin{\ast_{\kappa'}}}
\newcommand{\aststd}   {\mathbin{\ast_{\scriptscriptstyle\mathrm{Std}}}}
\newcommand{\astastd}  {\mathbin{\ast_{\cc{\scriptscriptstyle\mathrm{Std}}}}}
\newcommand{\astweyl}  {\mathbin{\ast_{\scriptscriptstyle\mathrm{Weyl}}}}
\newcommand{\fib}      {\scriptscriptstyle\mathrm{f\/ib}}
\newcommand{\fibweyl}  {\mathbin{\circ_{\scriptscriptstyle\mathrm{Weyl}}}}
\newcommand{\fibstd}   {\mathbin{\circ_{\scriptscriptstyle\mathrm{Std}}}}
\newcommand{\fibantistd}   {\mathbin{\circ_{\scriptscriptstyle\cc{\mathrm{Std}}}}}
\newcommand{\fibkap}   {\mathbin{\circ_\kappa}}
\newcommand{\fibkapStr}   {\mathbin{\circ_{\kappa'}}}

\newcommand{\adkap}    {\operatorname{\mathrm{ad}_\kappa}}
\newcommand{\adkapStr}    {\operatorname{\mathrm{ad}_{\kappa^\prime}}}
\newcommand{\taustd}   {\tau_{\scriptscriptstyle\mathrm{Std}}}
\newcommand {\W} [1] {\mathcal{W}_{#1}}
\newcommand {\WS} [2] {\mathcal{W}_{#1} \otimes \mathcal{S}^{#2}}
\newcommand {\WSL} [3] {\mathcal{W}_{#1} \otimes \mathcal{S}^{#2} \otimes \Lambda^{#3}}
\renewcommand {\S} [1] {\mathcal{S}^{#1}}
\newcommand {\FedDer} {{\mathcal{D}}}
\newcommand{\SymD}      {\mathsf{D}}


\newcommand{\refitem}[1]{$\ref{#1})$}
\numberwithin{equation}{section}

\begin{document}

\allowdisplaybreaks

\renewcommand{\thefootnote}{$\star$}

\renewcommand{\PaperNumber}{074}

\FirstPageHeading

\ShortArticleName{Deformation Quantization of Poisson Structures Associated to
  Lie Algebroids}

\ArticleName{Deformation Quantization of Poisson Structures\\ Associated to
  Lie Algebroids\footnote{This paper
is a contribution to the Special Issue on Deformation
Quantization. The full collection is available at
\href{http://www.emis.de/journals/SIGMA/Deformation_Quantization.html}{http://www.emis.de/journals/SIGMA/Deformation\_{}Quantization.html}}}

\Author{Nikolai NEUMAIER and Stefan WALDMANN}

\AuthorNameForHeading{N. Neumaier and S. Waldmann}

\Address{Fakult{\"a}t f{\"u}r Mathematik und Physik,
  Albert-Ludwigs-Universit{\"a}t Freiburg,\\
  Physikalisches Institut,
  Hermann Herder Stra{\ss}e 3,
  D-79104 Freiburg, Germany}

\Email{\href{mailto:Nikolai.Neumaier@physik.uni-freiburg.de}{Nikolai.Neumaier@physik.uni-freiburg.de},
\href{mailto:Stefan.Waldmann@physik.uni-freiburg.de}{Stefan.Waldmann@physik.uni-freiburg.de}}

\URLaddress{\url{http://idefix.physik.uni-freiburg.de/~nine/},\\
\hspace*{10.5mm}\url{http://idefix.physik.uni-freiburg.de/~stefan/}}

\ArticleDates{Received September 26, 2008, in f\/inal form May 25,
2009; Published online July 16, 2009}

\Abstract{In the present paper we explicitly construct deformation
    quantizations of certain Poisson structures on $E^*$, where $E
    \longrightarrow M$ is a Lie algebroid. Although the considered
    Poisson structures in general are far from being regular or even
    symplectic, our construction gets along without Kontsevich's
    formality theorem but is based on a generalized Fedosov
    construction. As the whole construction merely uses geometric
    structures of $E$ we also succeed in determining the dependence of
    the resulting star products on these data in f\/inding appropriate
    equivalence transformations between them. Finally, the
    concreteness of the construction allows to obtain explicit
    formulas even for a wide class of derivations and
    self-equivalences of the products.  Moreover, we can show that
    some of our products are in direct relation to the universal
    enveloping algebra associated to the Lie algebroid. Finally, we
    show that for a certain class of star products on $E^*$ the
    integration with respect to a density with vanishing modular
    vector f\/ield def\/ines a trace functional.}


\Classification{53D55; 53D17}

\section{Introduction}
\label{sec:Intro}

The question of existence and classif\/ication of formal star products
deforming an arbitrary Poisson bracket on a Poisson manifold has been
answered positively with the aid of Kontsevich's famous formality
theorem \cite{kontsevich:2003a}, see also
\cite{cattaneo.felder.tomassini:2002b} for a Fedosov-like
globalization. Besides the regular or even symplectic case, which was
very well understood long before, see \cite{dewilde.lecomte:1983b,
  fedosov:1994a, omori.maeda.yoshioka:1991a}, up to then the quantity
of examples of truly Poisson structures, that were known to be
deformation quantizable, was very poor: essentially the constant and
the linear Poisson brackets on a vector space. The linear Poisson
brackets are precisely the Kirillov--Kostant--Souriau brackets on the
dual $\mathfrak{g}^*$ of a Lie algebra $\mathfrak{g}$. A star product
for this Poisson bracket has been obtained in the early work of Gutt
\cite{gutt:1983a}.  Even after the general existence results the
number of examples that could be `handled' concretely remained small
whence there is a justif\/ied interest in such concrete examples.

In the present paper we consider Poisson brackets and their
deformation quantizations that in some sense generalize both canonical
Poisson brackets on cotangent bundles, see also
\cite{bordemann.neumaier.pflaum.waldmann:2003a,
bordemann.neumaier.waldmann:1999a,
bordemann.neumaier.waldmann:1998a}, and linear Poisson brackets on
$\mathfrak{g}^*$, see \cite{gutt:1983a}. These extreme examples can be
merged into the Poisson geometry of $(E^*,\theta_E)$, where $E
\longrightarrow M$ is a Lie algebroid and $\theta_E$ is the
corresponding Poisson tensor on the dual bundle. One main point of our
construction, that consists in a~modif\/ied Fedosov procedure, is that
it is \emph{elementary} in the sense that it only relies on tensor
calculus and avoids the use of the formality theorem. In addition,
there are even more Poisson brackets that depend on an additional
$\D_E$-closed $E$-two-form accessible by our method than just the
canonical one.  Moreover, it is geometric and permits to investigate
various properties of the constructed products.

Yet another independent motivation for this particular class of
Poisson structures and corresponding star products is given as
follows. Let $E^* \longrightarrow M$ be the dual of a vector bundle $E
\longrightarrow M$ and let $\xi \in \Gamma^\infty(TE^*)$ be the Euler
vector f\/ield on $E^*$.
\begin{Definition}[Homogeneous star product]
\label{definition:HomogeneousStarProduct}
    A star product $\star$ on $E^*$ is called homogeneous if
    $\mathsf{H} = \Lie_\xi + \lambda \frac{\partial}{\partial
      \lambda}$ is a derivation of $\star$.
\end{Definition}
Such homogeneous star products played a crucial role in the f\/irst
existence proofs on cotangent bundles, see \cite{cahen.gutt:1982a,
  dewilde.lecomte:1983a}. The homogeneous star products enjoy very
nice features which makes them an interesting class to study. The
following properties are verif\/ied in complete analogy to the case of
cotangent bundles:
\begin{Proposition}
    \label{proposition:HomogeneousStarProducts}
    Let $\star = \sum_{r=0}^\infty \lambda^r C_r$ be a homogeneous
    star product on $E^*$.
    \begin{enumerate}    \itemsep=0pt
    \item\label{item:HomogeneousStarLinearPoissonStructure} The
        Poisson structure $\{f, g\}_\star = \frac{1}{\I} (C_1(f, g) - C_1(g,
        f))$ is a linear Poisson structure whence $E$ is a Lie
        algebroid.
    \item\label{item:fstargPolynomials} Let $f, g \in
        \Pol^\bullet(E^*)$ be of degree $k$ and $\ell$, respectively.
        Then
        \begin{equation*}
            f \star g = \sum_{r=0}^{k+\ell} \lambda^r C_r(f, g)
            \qquad
            \textrm{with}
            \quad
            C_r(f, g) \in \Pol^{k+\ell-r}(E^*).
        \end{equation*}
    \item\label{item:PolynomialSubAlgebra} The polynomial functions
        $\Pol^\bullet(E^*)[\lambda]$ are a
        $\mathbb{C}[\lambda]$-algebra with respect to $\star$. In
        particular, $\star$ converges trivially on
        $\Pol^\bullet(E^*)[\lambda]$ for all $\lambda = \hbar \in
        \mathbb{C}$.
    \item\label{item:PolynomialGenerated} The polynomial functions
        $\Pol^\bullet(E^*)[\lambda]$ are generated by
        $\Pol^0(E^*)[\lambda]$ and $\Pol^1(E^*)[\lambda]$ with respect
        to $\star$.
    \item\label{item:PolNullCommutative} $\Pol^0(E^*)$ is a subalgebra
        with $\star$ being the undeformed commutative product.
    \end{enumerate}
\end{Proposition}
It was brought to our attention by Simone Gutt that the homogeneity
arguments from \cite{cahen.gutt:1982a, dewilde.lecomte:1983a} should
also apply to our more general situation beyond the cotangent bundle
case yielding an existence proof independent of our approach.

Following a suggestion of the referee, one may also wonder whether on
a particular geometry like a vector bundle $E^* \longrightarrow M$
there is a simplif\/ied construction of a formality which would allow
for a construction of a star product for the linear Poisson structure
analogously to the approach of Dito \cite{dito:1999a} or Kathotia
\cite{kathotia:2000a}. We leave this question for a future
investigation.

For the case of an integrable Lie algebroid there is yet another
approach to quantization namely that of \emph{strict} deformation
quantization considered in \cite{landsman.ramazan:2001a} which is
based on a generalized Weyl quantization using pseudo-dif\/ferential
operators. Also for the non-integrable case using the local
integrating Lie groupoids, the latter product has been related to the
universal enveloping algebra associated to the Lie algebroid in
\cite{nistor.weinstein.xu:1999a}. It yields a star product on the
polynomial functions on $E^*$. Actually, it is \emph{not} evident, see
e.g.~\cite{bordemann.neumaier.waldmann:1998a}, but in fact true, that
this product extends to a well-def\/ined bidif\/ferential formal star
product on \emph{all} smooth functions $C^\infty(E^*)$. Also Chemla
discussed, in the framework of complex Lie algebroids, the relation
between Lie algebroids and certain algebras of dif\/ferential operators,
see e.g.\ \cite{chemla:1999a}. Finally, one should note that although
there are some similarities our setting and our construction are
dif\/ferent from that in \cite{nest.tsygan:2001a}, where
\emph{symplectic} Lie algebroids are quantized and not the dual $E^*$
with its in general truly non-symplectic Poisson structure.

Our paper is organized as follows: in Section~\ref{sec:Preliminaries}
we collect some preliminaries on Lie algebroids $E \longrightarrow M$
and associated Poisson brackets on $E^*$ to f\/ix our notation. In
Section~\ref{sec:Fedosov} we present our construction of star products
on the \emph{polynomial} functions on $E^*$. We show that they
actually quantize the given Poisson brackets and that they extend in a
unique way to star products on all smooth functions. Our construction
depends on an ordering parameter $\kappa$, an $E$-connection $\nabla$,
and a formal series of $\D_E$-closed $E$-two-forms $B$.
Section~\ref{sec:Properties} is devoted to a further investigation of
the products. In particular, we show that among our star products
there are \emph{homogeneous} star products.  In
Section~\ref{sec:EquivalenceTransformations} we explicitly construct
equivalence transformations between the products provided the
parameters are appropriately related.  In
Section~\ref{sec:Representations} we make contact to the results of
\cite{nistor.weinstein.xu:1999a} and relate the products obtained in
the particular case $B=0$ to the universal enveloping algebra of $E$.
Finally, Section~\ref{sec:Trace} shows that in the case of a
\emph{unimodular} Lie algebroid the classical trace consisting in the
integration with respect to a certain constant density on $E^*$ is
also a trace functional for \emph{any} homogeneous star product on
$E^*$.  Here we establish a relation between the a priori
\emph{different} notions of unimodularity of the Lie algebroid $E$ and
the Poisson manifold $(E^*,\theta_E)$ that is interesting for its own,
see also \cite[Section~7]{weinstein:1997a}.

\section{Preliminaries on linear Poisson structures and Lie algebroids}
\label{sec:Preliminaries}

In this section we collect some well-known facts on linear Poisson
structures and Lie algebroids in order to f\/ix our notation. For
details see e.g.\ \cite{mackenzie:2005a, landsman:1998a,
  cannasdasilva.weinstein:1999a}.

Let $E \longrightarrow M$ be a vector bundle of f\/ibre dimension $N$
over an $n$-dimensional manifold. Then~$E$ is a \emph{Lie algebroid}
if it is equipped with a bundle map, the \emph{anchor} $\varrho: E
\longrightarrow TM$, and a~Lie bracket $[\cdot, \cdot]_E$ for the
sections $\Gamma^\infty(E)$ in $E$ such that we have the Leibniz rule
$[s, ut]_E = u [s, t]_E + (\varrho(s)u) t$ for all $u \in C^\infty(M)$
and $s, t \in \Gamma^\infty(E)$. It follows that $\varrho([s, t]_E) =
[\varrho(s), \varrho(t)]$. There are (among many others) three
examples of interest:
\begin{Example}[Lie algebroids]\qquad
    \label{example:LieAlgebroids}
    \begin{enumerate}
    \itemsep=0pt
    \item\label{item:TMLieAlgebroid} Clearly $TM$ with $\varrho = \id$
        and $[\cdot, \cdot]_{TM}$ the canonical Lie bracket is a Lie
        algebroid. Thus Lie algebroids over $M$ generalize the tangent
        bundle of $M$.
    \item\label{item:LieAlgebraLieAlgebroid} If $M = \{\mathrm{pt}\}$
        is a point then a Lie algebroid $\mathfrak{g} \longrightarrow
        \{\mathrm{pt}\}$ is nothing but a Lie algebra. Thus Lie
        algebroids generalize Lie algebras, also.
    \item\label{item:PoissonLieAlgebroid} If $(M, \theta)$ is a
        Poisson manifold then $T^*M$ becomes a Lie algebroid with
        anchor $\varrho(\alpha) = - \alpha^\varsharp$ and bracket
        $[\alpha, \beta]_\theta = - \Lie_{\alpha^\varsharp} \beta +
        \Lie_{\beta^\varsharp} \alpha - \D(\theta(\alpha, \beta))$,
        where $\alpha, \beta \in \Gamma^\infty(T^*M)$ and $^\varsharp:
        T^*M \longrightarrow TM$ is def\/ined by $\alpha^\varsharp =
        \theta(\cdot, \alpha)$.
    \end{enumerate}
\end{Example}

One main theme in the theory of Lie algebroids is to replace the
tangent bundle by $E$ and to translate geometric concepts based on the
(co-)tangent bundle into the language of Lie algebroids: Indeed, for a
Lie algebroid $E$, the dual bundle $E^*$ becomes a Poisson manifold
with a Poisson bracket $\{\cdot, \cdot\}_E$ such that
\begin{equation}
    \label{eq:PoissonPolynomial}
    \{\Pol^k(E^*), \Pol^\ell(E^*)\}_E \subseteq \Pol^{k+\ell-1}(E^*),
\end{equation}
where $\Pol^k(E^*) \subseteq C^\infty(E^*)$ denotes those functions on
$E^*$ which are homogeneous polynomials in f\/ibre direction of degree
$k$. This Poisson bracket is explicitly determined by
\begin{gather}
    \label{eq:LinearPoissonBracket}
    \{\pi^*u, \pi^*v\}_E = 0,
    \quad
    \{\pi^*u, \mathcal{J}(s)\}_E = \pi^*(\varrho(s)u),
    \quad
    \textrm{and}
    \quad
    \{\mathcal{J}(s), \mathcal{J}(t)\}_E = - \mathcal{J}([s, t]_E),
\end{gather}
where $\mathcal{J}: \mathcal{S}^\bullet(E) = \bigoplus_{k=0}^\infty
\Gamma^\infty(\Sym^k E) \longrightarrow \Pol^\bullet(E^*)$ denotes the
canonical graded algebra isomorphism between the symmetric $E$-tensor
f\/ields and polynomial functions on $E^*$. The signs in
\eqref{eq:LinearPoissonBracket} are convention and yield the `correct'
canonical Poisson structure of $T^*M$ in the case of $E = TM$, but the
negative of the usual Kirillov--Kostant--Souriau bracket on $\mathfrak{g}^*$
in the case of a Lie algebra $E = \mathfrak{g}$. This also motivates the
notion of \emph{linear} Poisson structures. The property
\eqref{eq:PoissonPolynomial} can equivalently be described in terms of
the corresponding Poisson tensor $\theta_E \in \Gamma^\infty(\Anti^2
TE^*)$ of $\{\cdot, \cdot\}_E$ and the Euler vector f\/ield $\xi \in
\Gamma^\infty(TE^*)$ on $E^*$ by
\begin{equation*}
    \Lie_\xi \theta_E = - \theta_E.
\end{equation*}

Moreover, for $\Gamma^\infty(\Anti^\bullet E^*)$ we obtain a
\emph{differential} $\D_E$, i.e.\ a super-derivation of the
$\wedge$-product of degree $+1$ with $\D_E^2 = 0$ by literally copying
the formula for the deRham dif\/ferential. In particular, for $u \in
C^\infty(M)$ we have $(\D_E u)(s) = \varrho(s) u$, and for $\alpha \in
\Gamma^\infty(E^*)$ one obtains
\begin{equation*}
    (\D_E\alpha)(s, t) = \varrho(s)(\alpha(t))  -
    \varrho(t)(\alpha(s)) - \alpha([s, t]_E),
\end{equation*}
where $s, t \in \Gamma^\infty(E)$. The corresponding cohomology theory
is the \emph{Lie algebroid cohomology}, denoted by
$\mathrm{H}^\bullet_E(M)$.  In fact, all these three structures, Lie
algebroid, linear Poisson structure, and dif\/ferential, are completely
equivalent, see e.g.\ \cite[Section~4.2]{waldmann:2006a:prebook} for a
pedagogical discussion.

\begin{Definition}[Gauged Poisson bracket]
    \label{definition:gauged}
    For an arbitrary $\D_E$-closed $E$-two-form $B_0$ we def\/ine the
    Poisson bracket $\{\cdot,\cdot\}_{B_0}$ by
\begin{equation*}
    \{\pi^*u, \pi^*v\}_{B_0} = \{\pi^*u, \pi^*v\}_E = 0,
    \qquad
    \{\pi^*u, \mathcal{J}(s)\}_{B_0} = \{\pi^*u, \mathcal{J}(s)\}_E
    = \pi^*(\varrho(s)u),
\end{equation*}
and
\begin{equation*}
    \{\mathcal{J}(s), \mathcal{J}(t)\}_{B_0}=
    \{\mathcal{J}(s), \mathcal{J}(t)\}_E  -
    \pi^* B_0(s,t)  = - \mathcal{J}([s, t]_E) -
    \pi^* B_0(s,t).
\end{equation*}
We call $\{\cdot,\cdot\}_{B_0}$ the \emph{gauged} Poisson bracket
(corresponding to $B_0$).
\end{Definition}
\begin{Remark}\qquad{}
    \begin{enumerate}
    \itemsep=0pt
    \item Obviously, $\{\cdot,\cdot\}_{B_0}$ is just a slight
        generalization of $\{\cdot,\cdot\}_E$ that coincides with this
        canonical Poisson bracket in case $B_0=0$. Observe that the
        $\D_E$-closedness of $B_0$ guarantees the Jacobi identity to
        be satisf\/ied.
    \item One should note that in general the Poisson bracket
        $\{\cdot,\cdot\}_{B_0}$ does \emph{not} coincide with the
        Poisson bracket obtained from $\theta_E$ via a gauge
        transformation $\Phi_F$ in the common sense (cf.~\cite{bursztyn.radko:2003a}) which arises from a $\D$-closed
        two-form $F$ on $E^*$. On the one hand there are Poisson
        brackets obtained from $\Phi_F(\theta_E)$ that are not of the
        particular form $\{\cdot,\cdot\}_{B_0}$ but on the other hand
        there are even Poisson brackets $\{\cdot,\cdot\}_{B_0}$ that
        cannot be obtained via such gauge transformations. For
        instance consider the Lie algebroid with $[\cdot,\cdot]_E =0$,
        then clearly $\theta_E =0$ and hence for all $F$ the Poisson
        tensor $\Phi_F(\theta_E)$ vanishes, but in contrast every
        $E$-two-form $B_0\ne 0$ is $\D_E$-closed and hence yields a
        non-vanishing Poisson bracket~$\{\cdot,\cdot\}_{B_0}$. Another
        in some sense opposite extreme case occurs in case $E = TM$
        with the canonical bracket, then it is easy to see that
        $\Phi_F (\theta_E)$ yields a bracket of the form
        $\{\cdot,\cdot\}_{B_0}$ if\/f~$F = \pi^*B_0$.
    \end{enumerate}
\end{Remark}

In the following we shall make use of covariant derivatives. The
adapted notion of a covariant derivative in the context of Lie
algebroids is an \emph{$E$-connection}, i.e.\ a bilinear map $\nabla:
\Gamma^\infty(E) \times \Gamma^\infty(E) \longrightarrow
\Gamma^\infty(E)$ such that $\nabla_s t$ is $C^\infty(M)$-linear in
$s$ and satisf\/ies the Leibniz rule $\nabla_s( u t) = (\varrho(s)u) t + u
\nabla_s t$, where $u \in C^\infty(M)$ and $s, t \in \Gamma^\infty(E)$.
For a detailed study of $E$-connections see \cite{fernandes:2002a}.
\begin{Example}
    \label{example:Econnection}
    Let $\nabla^E$ be a linear connection for $E$. Then $\nabla_s t =
    \nabla^E_{\varrho(s)} t$ is an $E$-connection.
\end{Example}

By compatibility with natural pairing and tensor products, $\nabla$
extends to all tensor powers of~$E$ and~$E^*$ in the usual way. On
functions we set $\nabla_s u = \varrho(s) u$. The \emph{torsion} $T
\in \Gamma^\infty(E \otimes \Anti^2 E^*)$ is def\/ined by
\begin{equation*}
    T(s, t) = \nabla_s t - \nabla_t s - [s, t]_E
\end{equation*}
and the \emph{curvature} $R \in \Gamma^\infty(\End(E) \otimes \Anti^2
E^*)$ is def\/ined by
\begin{equation*}
    R(s, t)
    = \nabla_s \nabla_t
    - \nabla_t \nabla_s
    - \nabla_{[s, t]_E}
\end{equation*}
as usual. If the torsion $T$ is non-trivial, one can pass to a new
$E$-connection by subtracting $\frac{1}{2}T$ to obtain a
\emph{torsion-free} $E$-connection. Note however, that if the
$E$-connection is of the form as in Example~\ref{example:Econnection},
then the result will in general no longer be of this form.

Finally, we sometimes make use of local formulas. By $e_1, \ldots, e_N
\in \Gamma^\infty(E\at{U})$ we denote a local basis of sections,
def\/ined on some suitable open subset $U \subseteq M$. Then $e^1,
\ldots, e^N \in \Gamma^\infty(E^*\at{U})$ denotes the corresponding
dual basis. We have induced linear f\/ibre coordinates $v^\alpha =
\mathcal{J}(e^\alpha) \in C^\infty(E \at{U})$ and $p_\alpha =
\mathcal{J}(e_\alpha) \in C^\infty(E^*\at{U})$.  If in addition $x^1,
\ldots, x^n$ are local coordinates on $U$ then the anchor and the
bracket are determined by
\begin{equation*}
    \varrho\at{U}
    = \varrho^i_\alpha e^\alpha
    \otimes \frac{\partial}{\partial x^i}
    \qquad
    \textrm{and}
    \qquad
    [e_\alpha, e_\beta]_E
    = c^\gamma_{\alpha\beta} e_\gamma
    \end{equation*}
with locally def\/ined functions $\varrho^i_\alpha,
c^\gamma_{\alpha\beta} \in C^\infty(U)$, respectively. Together with
the f\/ibre coordinates $p_1, \ldots, p_N$ the functions $q^1 = \pi^*
x^1, \ldots, q^n = \pi^*x^n$ provide a system of local coordinates on
$\pi^{-1}(U) \subseteq E^*$. Then the Poisson tensor is locally given
by
\begin{equation*}
    \theta_E \at{\pi^{-1}(U)} =
    \pi^* \varrho^i_\alpha \frac{\partial}{\partial q^i} \wedge
    \frac{\partial}{\partial p_\alpha}
    -
    \frac{1}{2} p_\gamma \pi^* c^\gamma_{\alpha\beta}
    \frac{\partial}{\partial p_\alpha}
    \wedge
    \frac{\partial}{\partial p_\beta}
\end{equation*}
and the \emph{Hamiltonian vector field} $X_f = \Schouten{f, \theta_E}$
of $f \in C^\infty(E^*)$ is
\begin{equation}
    \label{eq:HamiltonianVectorField}
    X_f \at{\pi^{-1}(U)}
    =
    - \pi^*\varrho^i_\alpha \frac{\partial f}{\partial q^i}
    \frac{\partial}{\partial p_\alpha}
    +
    \pi^*\varrho^i_\alpha \frac{\partial f}{\partial p_\alpha}
    \frac{\partial}{\partial q^i}
    +
    p_\gamma \pi^* c^\gamma_{\alpha\beta}
    \frac{\partial f}{\partial p_\alpha}
    \frac{\partial}{\partial p_\beta}.
\end{equation}

Given an $E$-connection $\nabla$ we locally have $\nabla_{e_\alpha}
e_\beta = \Gamma^\gamma_{\alpha\beta} e_\gamma$ with Christof\/fel
symbols $\Gamma^\gamma_{\alpha\beta} \in C^\infty(U)$. Then $\nabla$
is torsion-free if\/f locally $\Gamma_{\alpha\beta}^\gamma -
\Gamma_{\beta\alpha}^\gamma = c_{\alpha\beta}^\gamma$. For a section
$s \in \Gamma^\infty(E)$ we have a~\emph{horizontal lift} $s^\hor \in
\Gamma^\infty(TE^*)$ to $E^*$, locally given by
\begin{equation*}
    s^\hor \at{\pi^{-1}(U)} =
    \pi^*\left(s^\alpha \varrho_\alpha^i\right)
    \frac{\partial}{\partial q^i}
    +
    p_\gamma
    \pi^*\big(s^\alpha \Gamma_{\alpha\beta}^\gamma\big)
    \frac{\partial}{\partial p_\beta}.
\end{equation*}
Moreover, lifting horizontally is compatible with multiplication by
functions in the sense that $(u s)^\hor = (\pi^*u) s^\hor$.  The section
$\id \in \Gamma^\infty(\End(E)) \cong \Gamma^\infty(E \otimes E^*)$ can
now be lifted to a bivector f\/ield $\id^\lift \in
\Gamma^\infty(\bigotimes^2 TE^*)$ by lifting the $E^*$-part vertically
and the $E$-part horizontally, i.e.\ locally $\id^\lift =
(e_\alpha)^\hor \otimes (e^\alpha)^\ver$. For the linear Poisson
structure on $E^*$ we have
\begin{equation}
    \label{eq:ThetaEInTermsOfNabla}
    \theta_E = (e_\alpha)^\hor \wedge (e^\alpha)^\ver,
\end{equation}
i.e.\ the anti-symmetric part of $\id^\lift$ is $\theta_E$. Note that
$\theta_E$ does not depend on the choice of $\nabla$ but $\id^\lift$
does. Analogously, one has
\begin{equation}
    \label{eq:ThetaBInTermsOfNabla}
    \theta_{B_0} = \theta_E - B_0^\ver = (e_\alpha)^\hor \wedge
    (e^\alpha)^\ver - B_0^\ver
\end{equation}
for the Poisson tensor corresponding to the gauged Poisson bracket
$\{\cdot,\cdot\}_{B_0}$. This simple observation will be crucial
for the Fedosov construction of a star product on $E^*$.

%
%

\section{The Fedosov construction}
\label{sec:Fedosov}

The aim of this section is to obtain deformation quantizations for the
linear Poisson bracket on~$E^*$ and the Poisson brackets arising from
gauge transformations by $\D_E$-closed $E$-two forms in a geometric
fashion using a variant of Fedosov's construction. Since Fedosov's
construction is by now a very well-known approach to deformation
quantization, which has been adapted to many contexts as e.g.\
\cite{karabegov:2000a, bordemann.waldmann:1997a,
  bordemann.neumaier.waldmann:1999a,
  bordemann.neumaier.waldmann:1998a,neumaier:2003a}, we can be brief.
The version we are interested in resembles very much the approach of
\cite[Section~3]{bordemann.neumaier.waldmann:1999a} for cotangent
bundles.

Since a star product on $E^*$ is completely f\/ixed by its values on
$\Pol^\bullet(E^*)$ and since $\Pol^\bullet(E^*) \cong
\mathcal{S}^\bullet(E)$ via $\mathcal{J}$, we will construct a
deformation of the latter. We consider
\begin{equation}
    \label{eq:WSLDef}
    \mathcal{W}^\bullet \otimes \mathcal{S}^\bullet
    \otimes \Anti^\bullet
    =
    \prod_{k=0}^\infty \bigoplus_{\ell = 0}^\infty
    \Gamma^\infty\big(
        \Sym^k E^* \otimes \Sym^\ell E \otimes \Anti^\bullet  E^*
    \big)[[\lambda]],
\end{equation}
where we shall always use the sections of the complexif\/ied bundles in
view of applications in physics. Note however, that in the following
replacing the combination $\I \lambda$ by $\nu$ yields an entirely
real construction in the rescaled formal parameter $\nu$. Then
$\WSL{}{}{}$ is an associative $\mathbb{C}[[\lambda]]$-algebra with
respect to the symmetric/anti-symmetric tensor product $\mu$.

The three gradings in \eqref{eq:WSLDef} give rise to degree
derivations which we denote by $\degs$ for the $\Sym^\bullet
E^*$-degree, $\degsst$ for the $\Sym^\bullet E$-degree, and $\dega$
for the $\Anti^\bullet E^*$-degree, respectively. Moreover, we have
the $\lambda$-degree $\deg_\lambda = \lambda \frac{\partial}{\partial
  \lambda}$.  Using the symmetric and anti-symmetric insertion
derivations we locally have
\begin{gather*}
    \degs = (e^\alpha\! \otimes 1 \otimes 1) \inss(e_\alpha),
    \quad\!\!
    \degsst= (1 \otimes e_\alpha \otimes 1) \inss(e^\alpha),
    \quad\!\!
    \textrm{and}
    \quad\!\!
    \dega = (1 \otimes 1 \otimes e^\alpha) \insa(e_\alpha).
\end{gather*}
Moreover, we shall use the \emph{total degree} $\Deg = \degs +
\deg_\lambda$ as well as the \emph{homogeneity operator} (or:
$\lambda$-\emph{Euler derivation}) $\Eulerfib = \degsst +
\deg_\lambda$. The product $\mu$ of $\mathcal{W}^\bullet \otimes
\mathcal{S}^\bullet \otimes \Anti^\bullet$ is graded with respect to
all degree maps, i.e.\ $\degs$, $\dega$, $\degsst$, $\deg_\lambda$ and
hence $\Deg$ and $\Eulerfib$ are derivations. Moreover, $\mu$ is
super-commutative with respect to $\dega$.  The space of elements in
$\mathcal{W} \otimes \mathcal{S} \otimes \Anti$ of total degree $\geq
k$ will be denoted by $\WSL{k}{}{}$.

As usual in the Fedosov framework we consider the maps
\begin{equation*}
    \delta = (1\otimes 1 \otimes e^\alpha) \inss(e_\alpha)
    \qquad
    \textrm{and}
    \qquad
    \delta^* = (e^\alpha \otimes 1 \otimes 1) \insa(e_\alpha)
\end{equation*}
together with $\delta^{-1}$, def\/ined by $\delta^{-1}a =
\frac{1}{k+\ell} a$ for $\degs a = ka$ and $\dega a = \ell a$ if
$k+\ell > 0$ and $\delta^{-1}a = 0$ otherwise. Moreover,
\begin{equation*}
    \sigma: \ \mathcal{W} \otimes \mathcal{S} \otimes \Anti
    \longrightarrow \mathcal{S}
\end{equation*}
denotes the projection onto the part of $\degs$- and $\dega$-degree
$0$. Then one has
\begin{equation}
    \label{eq:DeltaHomotopies}
    \delta^2 = 0,
    \quad
    (\delta^{-1})^2 = 0
    \qquad
    \textrm{and}
    \qquad
    \delta \delta^{-1} + \delta^{-1} \delta + \sigma = \id.
\end{equation}
The degrees of the maps $\delta$ and $\delta^{-1}$ are
\begin{equation*}
    [\degs, \delta] = - \delta,
    \qquad
    [\dega, \delta] = \delta,
    \qquad
    \textrm{and}
    \qquad
    [\deg_\lambda, \delta] = 0 = [\degsst, \delta]
\end{equation*}
as well as
\begin{equation*}
    [\degs, \delta^{-1}] = \delta^{-1},
    \qquad
    [\dega, \delta^{-1}] = - \delta^{-1},
    \qquad
    \textrm{and}
    \qquad
    [\deg_\lambda, \delta^{-1}] = 0 = [\degsst, \delta^{-1}].
\end{equation*}
In a next step, one deforms $\mathcal{W} \otimes \mathcal{S} \otimes
\Anti$ using `$\kappa$-ordered' deformations, explicitly given by
\begin{equation*}
    a \fibkap b =
    \mu \circ \exp \left(
        - (1 -\kappa) \I\lambda \inss(e^\alpha) \otimes
        \inss(e_\alpha)
        + \kappa \I\lambda \inss(e_\alpha) \otimes \inss(e^\alpha)
    \right) (a \otimes b),
\end{equation*}
where $\kappa\in \mathbb{R}$ denotes a real parameter. Note that
$\fibkap$ is globally well def\/ined. We are mainly interested in
$\kappa \in [0,1]$.  Analogously to the case of cotangent bundles the
products for $\kappa = 0$, $\kappa =\frac{1}{2}$, and $\kappa=1$ are
referred to as standard-ordered, Weyl-ordered, and
anti-standard-ordered f\/ibrewise product and are denoted by
$\circ_0=\fibstd$, $\circ_{\frac{1}{2}}=\fibweyl$, and $\circ_1=
\fibantistd$, respectively. By the usual `commuting derivation'
argument, $\fibkap$ is associative for all $\kappa \in \mathbb{R}$.
All these deformed products quantize the same f\/ibrewise (super-)
Poisson structure $\{\cdot, \cdot\}_{\fib}$.

The total degree $\Deg$, the homogeneity operator $\Eulerfib$ as well
as $\dega$ are still derivations of $\fibkap$. Moreover,
\begin{equation*}
    \delta = \frac{\I}{\lambda} \adkap(1 \otimes \id_E),
\end{equation*}
where we view $\id_E \in \Gamma^\infty(\End(E))$ as an element in
$\mathcal{W}^0 \otimes \mathcal{S}^1 \otimes \Anti^1$ and $\adkap$
denotes the $\fibkap$-super-commutator. Thus $\delta$ is a quasi-inner
derivation and clearly $\delta (1 \otimes \id_E) = 0$.

Considering the f\/ibrewise `Laplacian' $\Delta_{\fib}$ def\/ined by
\begin{equation*}
    \Delta_{\fib} = \inss(e^\alpha)\inss(e_\alpha)
\end{equation*}
we have the identity
\begin{equation*}
    \Delta_{\fib} \circ \mu = \mu \circ \left(
        \Delta_{\fib} \otimes \id + \inss(e^\alpha)\otimes
        \inss(e_\alpha) + \inss(e_\alpha) \otimes \inss(e^\alpha)+ \id
        \otimes \Delta_{\fib}
    \right).
\end{equation*}
With this equation one can easily show that
\begin{equation}
    \label{eq:faserweiserNeumaier}
    \mathcal{M}_{\kappa' - \kappa} = \exp\left(
        \I \lambda (\kappa' -\kappa) \Delta_{\fib} \right)
\end{equation}
def\/ines a f\/ibrewise equivalence transformation from
$(\WSL{}{}{},\fibkap)$ to $(\WSL{}{}{},\fibkapStr)$. Note that we have
the commutation relations
\begin{equation*}
     [\delta ,\Delta_{\fib}] =
     [\deg_\lambda, \Delta_{\fib} ] =
     [\dega, \Delta_{\fib} ] = 0
     \qquad
     \textrm{and}
     \qquad
     [\degs, \Delta_{\fib} ]
     = [\degsst, \Delta_{\fib} ]
     = - \Delta_{\fib}.
\end{equation*}
The super-centre of $\fibkap$ is described as follows:
\begin{Lemma}
    \label{lemma:SuperCenter}
    $\adkap(a) = 0$ iff $\degs a = 0 = \degsst a$, i.e.\ $a \in
    \Gamma^\infty(\Anti^\bullet E^*)[[\lambda]]$.
\end{Lemma}

Now we choose a torsion-free $E$-connection $\nabla$ which gives an
exterior covariant derivative
\begin{equation*}
    D = (1 \otimes 1 \otimes e^\alpha) \nabla_{e_\alpha}.
\end{equation*}
Then $D$ is clearly globally well def\/ined. Moreover, and this is
remarkable compared to the usual Fedosov approach, $D$ is a
super-derivation of $\fibkap$ without any further assumptions. This is
clear, as $\fibkap$ only involves natural pairings.  The properties of
$D$ are easily determined and will be summarized in the following
lemma:
\begin{Lemma}
    \label{lemma:D}
    $D$ is a super-derivation of $\fibkap$ with $[\degs, D] =
    [\deg_\lambda, D] = [\degsst, D] = 0$, $[\dega, D]$ $=
    D$, and $[\Delta_{\fib},D]=0$. Moreover,
    \begin{equation*}
        D^2 = \frac{1}{2} [D, D] = - \frac{\I}{\lambda} \adkap(R)
        \qquad
        \textrm{and}
        \qquad
        [\delta, D] = 0.
    \end{equation*}
    The curvature $R \in \mathcal{W}^1 \otimes \mathcal{S}^1 \otimes
    \Anti^2$ satisfies the Bianchi identities $\delta R = 0 = DR$.
    Moreover, $\Delta_{\fib}R \in \Gamma^\infty (\Anti^2 E^*)$ is
    exact, i.e.\ there is an $E$-one-form $\Gamma$ such that $-
    \Delta_{\fib} R = \D_E \Gamma$. Finally, for $B \in
    \Gamma^\infty(\Anti^\bullet E^*)$ we have $D B = \D_E B$.
\end{Lemma}
Note that the lemma allows for some rather obvious modif\/ications in
case the $E$-connection has torsion analogous to
\cite{karabegov.schlichenmaier:2001b, neumaier:2001a}.

Now, one makes the usual ansatz
\begin{equation*}
    \FedDer_\kappa = - \delta + D + \frac{\I}{\lambda} \adkap
    (r_\kappa)
\end{equation*}
with an element $r_\kappa \in \WSL{1}{}{1}$. Clearly, $\FedDer_\kappa$
is a super-derivation of $\fibkap$, and we want to f\/ind~$r_\kappa$
such that $\FedDer_\kappa^2 = 0$. A direct computation yields that
\begin{equation*}
    \FedDer_\kappa^2 = \frac{\I}{\lambda}\adkap \left(- \delta
        r_\kappa + D r_\kappa + \frac{\I}{\lambda}
        r_\kappa \fibkap r_\kappa - R \right),
\end{equation*}
which vanishes if\/f $- \delta r_\kappa + D r_\kappa +
\frac{\I}{\lambda} r_\kappa \fibkap r_\kappa - R$ is a central element
in $\WSL{}{}{2}$. This is the case if\/f there is a formal series
of two-forms $B \in \Gamma^\infty(\Anti^2 E^*) [[\lambda]]$ with
\begin{equation*}
    \delta r_\kappa - D r_\kappa - \frac{\I}{\lambda} r_\kappa  \fibkap
    r_\kappa + R = B.
\end{equation*}
Since $[\FedDer_\kappa,[\FedDer_\kappa,\FedDer_\kappa]]=0$ by the
super-Jacobi identity the necessary condition for this equation to be
solvable is $\FedDer_\kappa B = \D_E B =0$.

After these preparations one is in the position to prove the following
theorem in analogy to \cite[Theorem~3.2]{fedosov:1994a} and
\cite[Theorem~5.2.2]{fedosov:1996a}:
\begin{Theorem}
    \label{theorem:rTHM}
    For every formal series $B= \sum_{j=0}^\infty \lambda^j B_j \in
    \Gamma^\infty(\Anti^2 E^*)[[\lambda]]$ of $\D_E$-closed
    $E$-two-forms there exists a unique element $r_\kappa \in
    \WSL{1}{}{1}$ such that
    \begin{equation}
        \label{eq:requations}
        \delta r_\kappa =  D r_\kappa +  \frac{\I}{\lambda}
        r_\kappa  \fibkap r_\kappa - R + B
        \qquad\textrm{and} \qquad
        \delta^{-1}r_\kappa =0.
    \end{equation}
    Moreover, $r_\kappa$ satisfies
    \begin{equation}
        \label{eq:rRecursion}
        r_\kappa = \delta^{-1}\left(D r_\kappa +
            \frac{\I}{\lambda} r_\kappa\fibkap r_\kappa  - R + B
        \right)
    \end{equation}
    from which $r_\kappa$ can be determined recursively. In this case
    the Fedosov derivation
    \begin{equation}
        \label{eq:FedDerivDef}
        \FedDer_\kappa = - \delta + D +
        \frac{\I}{\lambda} \adkap(r_\kappa)
    \end{equation}
    is a $\fibkap$-super-derivation of $\Anti^\bullet E^*$-degree $+1$
    and has square zero: $\FedDer_\kappa^2 = 0$.
\end{Theorem}
The proof follows the standard line of argument with the only
exception that thanks to our particular gradings also $B_0 \ne 0$ is
possible.

Before investigating the structure of $\ker(\FedDer_\kappa)\cap
\WS{}{}$ we note that the $\FedDer_\kappa$-cohomology is trivial on
elements $a$ with positive $\Anti^\bullet E^*$-degree since one has
the following homotopy formula
\begin{equation}
    \label{eq:DerHomotopy}
    \FedDer_\kappa \FedDer_\kappa^{-1} a + \FedDer_\kappa^{-1}
    \FedDer_\kappa a + \frac{1}{\id- [\delta^{-1}, D
      +\frac{\I}{\lambda}\adkap(r_\kappa)]} \sigma (a) = a,
\end{equation}
where
\begin{equation*}
    \FedDer_\kappa^{-1} a
    = -\delta^{-1}\left(
        \frac{1}{\id- [\delta^{-1}, D +
          \frac{\I}{\lambda}\adkap(r_\kappa)]}
        a \right)
\end{equation*}
for \emph{all} $a \in \WSL{}{}{}$ (cf.\
\cite[Theorem~5.2.5]{fedosov:1996a}). Clearly, \eqref{eq:DerHomotopy} is
a deformation of \eqref{eq:DeltaHomotopies}.

The next step in Fedosov's construction now consists in establishing a
bijection between the elements $a$ of $\WS{}{}$ with $\FedDer_\kappa a
=0$ and $\S{}$, which can be done as in \cite[Theorem~3.3]{fedosov:1994a}
or \cite[Theorem~5.2.4]{fedosov:1996a}.
\begin{Theorem}
    \label{theorem:FedosovTaylorTHM}
    Let $\FedDer_\kappa = - \delta + D +
    \frac{\I}{\lambda}\adkap(r_\kappa): \WSL{}{}{\bullet}
    \longrightarrow \WSL{}{}{\bullet+1}$ be given as in~\eqref{eq:FedDerivDef} with $r_\kappa$ as in~\eqref{eq:requations}.
    \begin{enumerate}
    \itemsep=0pt
    \item Then for any $s \in \S{} = \S{\bullet}(E)[[\lambda]]=
        \bigoplus_{\ell=0}^\infty \Gamma^\infty(\Sym^\ell
        E)[[\lambda]]$ there exists a unique element $\tau_\kappa(s)
        \in \ker (\FedDer_\kappa) \cap \WS{}{}$ such that
        \begin{equation*}
            \sigma (\tau_\kappa(s))= s
        \end{equation*}
        and $\tau_\kappa: \S{} \longrightarrow \ker (\FedDer_\kappa)
        \cap \WS{}{}$ is $\mathbb C[[\lambda]]$-linear and referred to
        as the Fedosov--Taylor series corresponding to $\FedDer_\kappa$.
    \item In addition $\tau_\kappa(s)$ is given by
        \begin{equation*}
            \tau_\kappa(s)= \frac{1}{\id- [\delta^{-1}, D +
              \frac{\I}{\lambda}\adkap(r_\kappa)]} s
        \end{equation*}
        or can be determined recursively from
        \begin{equation}
            \label{eq:GenTaylorRecurs}
            \tau_\kappa(s) = s + \delta^{-1}\left( D \tau_\kappa(s) +
                \frac{\I}{\lambda} \adkap(r_\kappa)\tau_\kappa (s)\right).
        \end{equation}
    \item Since $\FedDer_\kappa$ as constructed above is a
        $\fibkap$-super-derivation and since $\ker(\FedDer_\kappa) \cap
        \WS{}{}$ is a~$\fibkap$-subalgebra, a new associative product
        $\astkap$ for $\S{}$ that deforms the symmetric product
        is defined by pull-back of $\fibkap$ via $\tau_\kappa$:
        \begin{equation*}
            s \astkap t = \sigma ( \tau_\kappa(s)\fibkap \tau_\kappa(t)).
        \end{equation*}
    \end{enumerate}
\end{Theorem}

In the following we shall refer to the associative product $\astkap$
def\/ined above as the $\kappa$-ordered $E$-Fedosov product
(corresponding to $(\nabla, B)$).

Clearly, from the associative deformation above
$\mathcal{S}^\bullet(E)$ inherits the structure of a Poisson algebra
by the term occurring in the f\/irst order $\I \lambda$ of the
$\astkap$-commutator. As a f\/irst result we prove that the
corresponding Poisson bracket $\{\cdot,\cdot\}_{\astkap}$ coincides
with the one induced on $\mathcal{S}^\bullet(E)$ by pull-back of
$\{\cdot,\cdot\}_{B_0}$ via $\mathcal{J}$.
\begin{Proposition}
    \label{proposition:CorrectPoissonBracket}
    For all $s,t \in \mathcal{S}^\bullet(E)$ we have
    \begin{gather*}
        \! s \astkap t = s t
        + \I\lambda \!\left(\!
            - (1-\kappa) \inss (e^\alpha)s  \nabla_{e_\alpha}t\!
            + \kappa \nabla_{e_\alpha}s  \inss(e^\alpha)t\!
            - \tfrac{1}{2} B_0(e_\alpha,e_\beta)
            \inss(e^\alpha)s  \inss(e^\beta)t
        \right)\! + O(\lambda^2).\!
    \end{gather*}
    In particular, this implies that the induced Poisson bracket on
    $\mathcal{S}^\bullet(E)$ is given by
    \begin{equation*}
        \{s, t\}_{\astkap} = \nabla_{e_\alpha} s  \inss (e^\alpha) t -
        \inss (e^\alpha) s \nabla_{e_\alpha} t -
        B_0(e_\alpha,e_\beta) \inss(e^\alpha)s \inss(e^\beta)t =
        \mathcal{J}^{-1}\{\mathcal J (s),\mathcal J (t)\}_{B_0}.
    \end{equation*}
    In case $B_0= 0$ we therefore have $\{s,t\}_{\astkap}=
    \mathcal{J}^{-1}\{\mathcal{J} (s),\mathcal{J} (t)\}_E$.
\end{Proposition}
\begin{proof}
    By straightforward computation using that $\sigma$ commutes with
    $\inss (\alpha)$ for all $\alpha \in \Gamma^\infty(E^*)$ one f\/inds
    \[
    s \astkap t = s  t
    - (1 - \kappa) \I\lambda
    \inss (e^\alpha)s \inss(e_\alpha) \tau_\kappa (t)^1_\cl
    + \kappa \I \lambda
    \inss(e_\alpha) \tau_\kappa (s)^1_\cl  \inss (e^\alpha)t
     + O(\lambda^2),
    \]
    where $\tau_\kappa (t)^1_\cl$ denotes the term of $\Sym^\bullet E^*$-degree
    $1$ of the classical part of $\tau_\kappa(t)$. From equation~\eqref{eq:GenTaylorRecurs} one f\/inds that $\tau_\kappa(t)^1_\cl$
    satisf\/ies
    \[
    \tau_\kappa(t)^1_\cl = \delta^{-1}\left( Dt + \inss
    (e^\alpha)r_{\kappa, \cl}^0 \inss (e_\alpha)\tau_\kappa (t)^1_\cl -
    \inss (e^\alpha) t \inss (e_\alpha) r_{\kappa, \cl}^1 \right),
    \]
    where again $_\cl$ refers to the classical part and the upper
    index indicates the $\Sym^\bullet E^*$-degree. But from
    \eqref{eq:rRecursion} it is evident that $r_{\kappa, \cl}^0 = 0$
    which in turn implies that $r_{\kappa, \cl}^1 = \delta^{-1}B_0$.
    Therefore we explicitly get $\tau(t)^1_\cl = \delta^{-1} \left( D
        t -\frac{1}{2} \inss(e^\beta) t \insa (e_\beta) B_0\right) =
    (1 \otimes e^\alpha) \left( \nabla_{e_\alpha}t + \frac{1}{2}
        B_0(e_\alpha,e_\beta) \inss(e^\beta) t\right)$. Together with
    the above expression for $s \astkap t$ this implies the f\/irst
    statement of the proposition. For the proof of the formula for
    $\{s,t\}_{\astkap}$ it suf\/f\/ices to check the statement on the
    generators of~$\mathcal{S}^\bullet(E)$, i.e.\ on functions $u,
    v\in C^\infty(M)$ and on $s,t \in \Gamma^\infty(E)$. But this is
    straightforward using the properties of a torsion-free
    $E$-connection. This fact can also directly be read of\/f the
    formulas~\eqref{eq:ThetaEInTermsOfNabla} and~\eqref{eq:ThetaBInTermsOfNabla}.
\end{proof}

Clearly, the products $\astkap$ give rise to star products $\starkap$
on $\Pol^\bullet(E^*)[[\lambda]]$ by def\/ining
\begin{equation*}
    f \starkap g = \mathcal{J}(\mathcal{J}^{-1}(f) \astkap \mathcal{J}^{-1}(g))
\end{equation*}
for $f,g\in \Pol^\bullet(E^*)[[\lambda]]$. According to
Proposition~\ref{proposition:CorrectPoissonBracket} this is a
deformation in the direction of the Poisson bracket
$\{\cdot,\cdot\}_{B_0}$. But in order to be able to extend the latter
product to a star product on $C^\infty(E^*)[[\lambda]]$ one
additionally has to show that on $\Pol^\bullet(E^*)[[\lambda]]$ the
product $\starkap$ can be described by \emph{bidifferential} operators
(cf.\ the discussion in
\cite[Section~3]{bordemann.neumaier.waldmann:1999a}).
\begin{Proposition}
    \label{proposition:TauDifferOrd}
    Writing the term of total degree $k \in \mathbb{N}$ in
    $\tau_\kappa(t)$ for $t \in \S{\bullet}(E)= \Gamma^\infty
    (\Sym^\bullet E)$ as $\tau_\kappa(t)^{(k)}= \sum_{\ell=0}^k
    \left(\frac{\lambda}{\I}\right)^{k-\ell} \tau_\kappa(t)^{(k),
      \ell}$ with $\tau_\kappa(t)^{(k), \ell} \in \Gamma^\infty\left(
        \Sym^\ell E^* \otimes \Sym^\bullet E \right)$ one has that the
    mapping
    \begin{equation*}
        t \mapsto  \tau_\kappa(t)^{(k),\ell}
    \end{equation*}
    is a differential operator of order $k$ for all $0\leq \ell \leq
    k$.
\end{Proposition}
\begin{proof}
    For the proof we use the usual algebraic def\/inition of
    dif\/ferential operators on appropriate subspaces of the
    $\S{\bullet}(E)= \Gamma^\infty(\Sym^\bullet E)$ (left) module
    $\prod_{k=0}^\infty \Gamma^\infty\left( \Sym^k E^* \otimes
        \Sym^\bullet E \otimes \Anti^\bullet E^* \right)$ that take
    their values again in some subspace of $\prod_{k=0}^\infty
    \Gamma^\infty\left( \Sym^k E^* \otimes \Sym^\bullet E \otimes
        \Anti^\bullet E^* \right)$, see e.g.\
    \cite{waldmann:2006a:prebook}.  Clearly, this def\/inition also
    applies to dif\/ferential operators on the algebra $\S{\bullet}(E)$
    with values in some subspace of the $\S{\bullet}(E)$ (left) module
    $\prod_{k=0}^\infty \Gamma^\infty\left( \Sym^k E^* \otimes
        \Sym^\bullet E \otimes \Anti^\bullet E^* \right)$. With this
    notion it is evident that $\delta^{-1}$ is a dif\/ferential operator
    of order $0$ since it commutes with all (left) multiplications by
    elements of $\S{\bullet}(E)$. Likewise, $\inss(s)$ is a
    dif\/ferential operator of order $0$ for all $s\in
    \Gamma^\infty(E)$. Furthermore, $D$ as well as $\inss(\alpha)$ for
    $\alpha \in \Gamma^\infty(E^*)$ are dif\/ferential operators of
    order $1$ since the respective commutators with (left)
    multiplications by elements of $\S{\bullet}(E)$ are given by
    (left) multiplications, which are by def\/inition dif\/ferential
    operators of order~$0$. Moreover, we use the well-known fact that
    composing two dif\/ferential operators of order $k$ and $k'$ the
    result is a dif\/ferential operator of order $k+k'$. After these
    preparations the proof of the proposition is a rather
    straightforward induction on the total degree.  It is trivial that
    $t \mapsto \tau_\kappa(t)^{(0)} = t$ is a dif\/ferential operator of
    order~$0$. Now let us assume that we have shown that for $j =
    0,\ldots,k$ we have that $t\mapsto \tau_\kappa(t)^{(j)}$ consists
    of dif\/ferential operators of order $j$. Using the recursion
    for\-mu\-la~\eqref{eq:GenTaylorRecurs} for $\tau_\kappa(t)$ and
    observing the respective total degree of the involved elements it
    is lengthy but straightforward to show that
    \begin{gather*}
            \tau_\kappa(t)^{(k+1)} = \delta^{-1}
            D \tau_\kappa(t)^{(k)}\\
            + \delta^{-1}\Bigg(
            \sum_{\ell=0}^k\sum_{j=0}^{k-\ell}\sum_{m=1}^{\ell+1}
            \left(\frac{\lambda}{\I}\right)^{j+m-1}
            \frac{(1-\kappa)^{j} (-\kappa)^m}{j!m!}
            \inss(e^{\alpha_1})\cdots \inss(e^{\alpha_{j}})
            \inss(e_{\beta_1})\cdots\inss(e_{\beta_m})
            r_\kappa^{(\ell+1)}\\
            \qquad{} \inss(e_{\alpha_1})\cdots\inss(e_{\alpha_{j}})
            \inss(e^{\beta_1})\cdots \inss(e^{\beta_m})
            \tau_\kappa(t)^{(k-\ell)}
            \Bigg)\\
            - \delta^{-1}\Bigg(
            \sum_{\ell=0}^k\sum_{j=1}^{\ell+1}\sum_{m=0}^{k-\ell}
            \left(\frac{\lambda}{\I}\right)^{j+m-1}
            \frac{(1-\kappa)^{j} (-\kappa)^m}{j!m!}
            \inss(e^{\alpha_1})\cdots \inss(e^{\alpha_{j}})
            \inss(e_{\beta_1})\cdots\inss(e_{\beta_m})
            \tau_\kappa(t)^{(k-\ell)}\\
           \qquad{} \inss(e_{\alpha_1})\cdots\inss(e_{\alpha_{j}})
            \inss(e^{\beta_1})\cdots \inss(e^{\beta_m})
            r_\kappa^{(\ell+1)} \Bigg).
        \end{gather*}
    Since $\delta^{-1}$ is a dif\/ferential operator of order $0$ and
    $D$ is a dif\/ferential operator of order $1$ the induction
    hypotheses yields that $t \mapsto \delta^{-1} D
    \tau_\kappa(t)^{(k)}$ only contains dif\/ferential operators of
    order at most $k+1$. Furthermore, from the fact that $t \mapsto
    \inss(e^{\beta_1})\cdots \inss(e^{\beta_m})
    \tau_\kappa(t)^{(k-\ell)}$ only contains dif\/ferential operators of
    order at most $k-\ell +m$ and observing that the sum over $m$ in
    the second summand runs from $1$ to $\ell+1$ the order of the
    dif\/ferential operators occurring in the second term is at most
    $k+1$.  Here we additionally have used that multiplications with
    elements in $\prod_{k=0}^\infty \Gamma^\infty\left( \Sym^k E^*
        \otimes \Sym^\bullet E \otimes \Anti^\bullet E^* \right)$ are
    of order $0$ and that the symmetric insertions of sections in $E$
    are of order $0$ also. The same line of argument applied to the
    third summand in $\tau_\kappa(t)^{(k+1)}$ yields that this part of
    the Fedosov Taylor series only consists in dif\/ferential operators
    of order $k+1$ proving the assertion by induction.
\end{proof}

After these preparations we can state the main result of this section:
\begin{Theorem}
    For all $\kappa \in \mathbb{R}$ and all $\D_E$-closed $B \in
    \Gamma^\infty(\Anti^2 E^*)[[\lambda]]$ one has:
    \begin{enumerate}
    \itemsep=0pt
    \item\label{item:astnatural} For all $s,t \in \S{\bullet}(E)$ the
        product $\astkap$ can be written as
        \begin{equation*}
            s \astkap t = \sum_{k=0}^\infty \lambda^k D_k(s,t),
        \end{equation*}
        where $D_k :\S{\bullet} (E) \times \S{\bullet}(E)
        \longrightarrow \S{\bullet}(E)$ is a bidifferential operator
        of order $k$ in both arguments.
    \item\label{item:starnatural} For all $f,g \in \Pol^\bullet(E^*)$
        the product $\starkap$ is given by
        \begin{equation*}
            f \starkap g = \sum_{k=0}^\infty \lambda^k C_k(f,g)
            = \sum_{k=0}^\infty \lambda^k
            \mathcal{J}(D_k(\mathcal{J}^{-1}(f),
            \mathcal{J}^{-1}(g))),
        \end{equation*}
        where $C_k :\Pol^\bullet(E^*) \times \Pol^\bullet(E^*)
        \longrightarrow \Pol^\bullet(E^*)$ is a bidifferential
        operator of order $k$ in both arguments. Hence $\starkap$
        extends in a unique way to a star product on
        $C^\infty(E^*)[[\lambda]]$ with respect to the Poisson bracket
        $\{\cdot,\cdot\}_{B_0}$ which is in addition natural in the
        sense of {\rm \cite{gutt.rawnsley:2003a}}.
    \end{enumerate}
\end{Theorem}
\begin{proof}
    For the proof of the f\/irst statement an easy computation
    expressing $s \astkap t$ by means of $\tau_\kappa(s)^{(m),\ell}$ and
    $\tau_\kappa(t)^{(k-m),\ell'}$ as in
    Proposition~\ref{proposition:TauDifferOrd} yields
    \begin{gather*}
            D_k(s,t) = (- \I)^k \sum_{m=0}^k \sum_{\ell'=0}^{k-m}
            \sum_{\ell=0}^m \frac{(1- \kappa)^{\ell'} (-\kappa)^\ell}{\ell!
            \ell'!}
             \inss(e^{\alpha_1}) \cdots \inss(e^{\alpha_{\ell'}})
            \inss(e_{\beta_1})\cdots \inss (e_{\beta_\ell})
            \tau_\kappa(s)^{(m),\ell}\\
           \phantom{D_k(s,t) = (- \I)^k}{} \inss(e^{\beta_1})\cdots \inss (e^{\beta_\ell})
            \inss(e_{\alpha_1}) \cdots \inss(e_{\alpha_{\ell'}})
            \tau_\kappa(t)^{(k-m),\ell'}.
        \end{gather*}
    Now according to Proposition~\ref{proposition:TauDifferOrd} the
    map $s \mapsto \inss(e^{\alpha_1}) \cdots
    \inss(e^{\alpha_{\ell'}}) \tau_\kappa(s)^{(m),\ell}$ is a
    dif\/ferential ope\-ra\-tor of order $m+\ell'$ and $t \mapsto
    \inss(e^{\beta_1})\cdots \inss (e^{\beta_\ell})
    \tau_\kappa(t)^{(k-m),\ell'}$ is a dif\/ferential operator of order
    $k-m+\ell$. Observing that the sum over $\ell'$ runs from $0$ to
    $k - m$ and that the sum over $\ell$ runs from $0$ to $m$ it is
    then obvious that the highest order of dif\/ferentiation occurring
    in $D_k(s,t)$ is $k$ in both arguments proving
    \refitem{item:astnatural}. The f\/irst part of the second statement
    is trivial since $\mathcal{J}$ is an isomorphism between the
    undeformed associative algebras $\S{\bullet}(E)$ and
    $\Pol^\bullet(E^*)$. The last part of the second statement follows
    from the fact that bidif\/ferential operators on $C^\infty(E^*)$ are
    completely determined by their values on $\Pol^\bullet(E^*)$.
\end{proof}
\begin{Remark}
    Because of the statements of the preceding theorem it is
    suf\/f\/icient to prove properties of $\starkap$ on
    $\Pol^\bullet(E^*)$ in order to show that they actually hold on
    all of $C^\infty(E^*)$. This observation drastically simplif\/ies
    many of the further investigations in the following sections.
\end{Remark}

\section[Further properties of $\starkap$]{Further properties of $\boldsymbol{\starkap}$}
\label{sec:Properties}

In this section we shall f\/irst f\/ind necessary conditions on the data
determining the star products~$\starkap$ that guarantee that these
products are homogeneous in the sense of
Def\/inition~\ref{definition:HomogeneousStarProduct}. These star
products are of particular interest since the results of
Section~\ref{sec:Trace} about the trace functional apply. Moreover, we
will show that the recursion formula for the element $r_\kappa$
determining the Fedosov derivation drastically simplif\/ies due to the
special shape of the f\/ibrewise product $\fibkap$ and provides some
rather explicit formulas. In particular, we explicitly compute the
star products of an arbitrary function $f \in C^\infty(E^*)$ with the
pull-back $\pi^*u$ of a function $u \in C^\infty(M)$ proving that
among our star products $\starkap$ there are star products of
(anti-)standard-ordered type.  Finally, we f\/ind conditions on which
the constructed star products are of Weyl type.
\begin{Proposition}
    \label{proposition:HomogenesStarKap}
    Let $\kappa \in \mathbb{R}$.
    \begin{enumerate}
    \itemsep=0pt
    \item\label{item:homogenBedingunganB} The element $r_\kappa \in
        \WSL{1}{}{1}$ constructed in Theorem~{\rm \ref{theorem:rTHM}}
        satisfies $\Eulerfib r_\kappa = r_\kappa$ iff $\deg_\lambda B
        = B$, which is the case iff $B = \lambda B_1$.  In this case,
        $\Eulerfib$ commutes with $\FedDer_\kappa$ implying that
        $\Eulerfib \tau_\kappa(s) = \tau_\kappa(\Eulerfib s)$ for all
        $s \in \S{}$ which in turn shows that
        \begin{equation}
            \label{eq:AstKapHomogen}
            \Eulerfib (s \astkap t ) = (\Eulerfib s) \astkap t + s \astkap
            (\Eulerfib t)
        \end{equation}
        holds for all $s,t \in \S{}$.
    \item\label{item:starhomogen} In case $B$ satisfies $\deg_\lambda
        B = B$ the star product $\starkap$ on
        $C^\infty(E^*)[[\lambda]]$ obtained from the Fedosov
        construction is homogeneous, i.e.\
        \begin{equation}
            \label{eq:StarKapHomogen}
            \Euler (f \starkap g ) =  (\Euler f) \starkap g + f \starkap
            (\Euler g)
        \end{equation}
        for all $f,g\in C^\infty(E^*)[[\lambda]]$.
    \end{enumerate}
\end{Proposition}
\begin{proof}
    Part~\refitem{item:starhomogen} is obvious from the observation
    that $\Euler = \mathcal{J} \circ \Eulerfib \circ \mathcal{J}^{-1}$
    and the very def\/inition of $\starkap$. One just has to observe
    that the identity \eqref{eq:StarKapHomogen} is an equation between
    (bi)dif\/ferential operators which is satisf\/ied in case it is
    satisf\/ied on all of $\Pol^\bullet(E^*)$ which is the case
    according to part~\refitem{item:homogenBedingunganB}. Thus we are
    left with the proof of the f\/irst part of the proposition. Applying
    $\Eulerfib$ to the equations determining $r_\kappa$ and combining
    the resulting equations with the original ones one gets
    $\delta^{-1}(\Eulerfib r_\kappa - r_\kappa)=0$ and
    $\delta(\Eulerfib r_\kappa - r_\kappa) = D (\Eulerfib r_\kappa -
    r_\kappa)+ \frac{\I}{\lambda}\adkap(r_\kappa)(\Eulerfib r_\kappa -
    r_\kappa) + \deg_\lambda B - B$. Therefore $\Eulerfib r_\kappa =
    r_\kappa$ evidently implies $\deg_\lambda B = B$. Vice versa
    supposing $\deg_\lambda B = B$ the above identities combine to
    $\Eulerfib r_\kappa - r_\kappa = \delta^{-1}(D (\Eulerfib r_\kappa
    - r_\kappa)+ \frac{\I}{\lambda}\adkap(r_\kappa)(\Eulerfib r_\kappa
    - r_\kappa))$ which is a f\/ixed point equation with a unique
    solution that trivially is solved by zero and hence $\Eulerfib
    r_\kappa = r_\kappa$. Using that $\Eulerfib$ is a derivation of
    $\fibkap$ it is straightforward to verify that in this case
    $[\Eulerfib, \FedDer_\kappa]=0$. Applying this identity to
    $\tau_\kappa(s)$ one f\/inds $\FedDer_\kappa(\Eulerfib
    \tau_\kappa(s))=0$ and therefore using
    Theorem~\ref{theorem:FedosovTaylorTHM} we have $\Eulerfib
    \tau_\kappa(s) = \tau_\kappa (\sigma (\Eulerfib\tau_\kappa (s) ))
    = \tau_\kappa(\Eulerfib s )$, where the last equality follows from
    the fact that $\sigma$ commutes with $\Eulerfib$. Finally,
    equation \eqref{eq:AstKapHomogen} is a direct consequence of the
    very def\/inition of $\astkap$ and the identity for $\Eulerfib
    \tau_\kappa(s)$ just shown.
\end{proof}

Now we turn to the more detailed consideration of the recursion
formulas for $r_\kappa$.
\begin{Proposition}
    \label{proposition:rproperties}
    Let $\kappa \in \mathbb{R}$.
    \begin{enumerate}
    \itemsep=0pt
    \item\label{item:rsimpler} For all $\D_E$-closed $B \in
        \Gamma^\infty(\Anti^2 E^*)[[\lambda]]$ the element $r_\kappa
        \in \WSL{1}{}{1}$ is of $\Sym^\bullet E$-degree at most $1$
        and is determined by the simplified recursion formula
        \begin{equation}
            \label{eq:simpler}
            r_\kappa  =  \delta^{-1}\left(
                D r_\kappa + \inss(e^\alpha) r_\kappa \inss(e_\alpha)
                r_\kappa - R + B\right).
        \end{equation}
        In particular, this shows that $r_\kappa$ is actually
        independent of $\kappa$. Henceforth we thus can neglect the
        additional index $\kappa$ and write $r$ for $r_\kappa$.
    \item\label{item:rzerlegt} For all $B$ as above $r$ can be written
        as $r = \mathsf{r}_0 + \mathsf{r}_1$. Here $\mathsf{r}_1$
        denotes the solution of \eqref{eq:requations} for $B=0$, which
        is homogeneous of $\Sym^\bullet E$-degree $1$ and
        $\mathsf{r}_0$, which is homogeneous of $\Sym^\bullet
        E$-degree~$0$ and explicitly given by
        \begin{equation}
            \label{eq:rNull}
            \mathsf{r}_0 =
            \frac{1}{\id - \delta^{-1}
              \left(
                  D + \inss(e^\alpha) \mathsf{r}_1 \inss(e_\alpha)
              \right)
            }
            \delta^{-1} B
            =
            \frac{1}{\id -
              \left[\delta^{-1},
                  D + \inss(e^\alpha) \mathsf{r}_1 \inss(e_\alpha)
              \right]
            }
            \delta^{-1} B.
        \end{equation}
        In particular, $\mathsf{r}_0$ depends linearly on $B$.
    \end{enumerate}
\end{Proposition}
\begin{proof}
    We prove assertion~\refitem{item:rsimpler} by induction on the
    total degree.  From the recursion formula~\eqref{eq:rRecursion} we
    f\/ind $r_\kappa^{(1)}= \delta^{-1}B_0$, $r_\kappa^{(2)} =
    \delta^{-1}( D r_\kappa^{(1)} + \frac{\I}{\lambda} r_\kappa^{(1)}
    \fibkap r_\kappa^{(1)} - R + \lambda B_1)$ and
    \[
    r_\kappa^{(k+1)} = \delta^{-1}\left( D r_\kappa^{(k)} +
        \frac{\I}{\lambda} \sum_{\ell=1}^k r_\kappa^{(\ell)} \fibkap
        r_\kappa^{(k+1-\ell)} + \lambda^k B_k \right)
    \]
    for $k \geq 2$. Obviously, $r_\kappa^{(1)}$ is of $\Sym^\bullet
    E$-degree $0$. Using this observation together with the explicit
    shape for $\fibkap$ yields $r_\kappa^{(2)} = \delta^{-1}( D
    r_\kappa^{(1)} - R + \lambda B_1)$ which is easily seen to be of
    $\Sym^\bullet E$-degree at most~$1$. Now assume that we already
    have shown that $r_\kappa^{(j)}$ is of $\Sym^\bullet E$-degree at
    most $1$ for $j = 1,\ldots,k$, then again use of the explicit
    shape of $\fibkap$ yields by straightforward computation that
    \[
    r_\kappa^{(k+1)} = \delta^{-1}\left( D r_\kappa^{(k)} +
        \sum_{\ell=1}^k \inss(e^\alpha) r_\kappa^{(\ell)} \inss (e_\alpha
        )r_\kappa^{(k+1-\ell)} + \lambda^k B_k \right).
    \]
    But this expression is obviously of $\Sym^\bullet E$-degree at
    most $1$ in case $r_\kappa^{(1)},\ldots,r_\kappa^{(k)}$ are.
    Recollecting the terms of each total degree into one element
    $r_\kappa$ f\/inally shows the simplif\/ied recursion formula~\eqref{eq:simpler}. According to \refitem{item:rsimpler} we can
    write $r = \mathsf{r}_0 + \mathsf{r}_1$, where $\mathsf{r}_0$ and
    $\mathsf{r}_1$ are homogeneous of $\Sym^\bullet E$-degree~$0$ and~$1$. Inserting this decomposition into the above recursion formula
    yields considering the parts of degree~$1$ and~$0$ separately
    $\mathsf{r}_1 = \delta^{-1}(D \mathsf{r}_1 +
    \inss(e^\alpha)\mathsf{r}_1 \inss(e_\alpha)\mathsf{r}_1 -R)$,
    which coincides with the equation for $r$ in case $B=0$, and
    $\mathsf{r}_0 = \delta^{-1}(D \mathsf{r}_0 + \inss(e^\alpha)
    \mathsf{r}_1 \inss (e_\alpha) \mathsf{r}_0 + B)$. This equation is
    uniquely solved by \eqref{eq:rNull}.
\end{proof}
\begin{Remark}
    One should note, that the homogeneity of $\mathsf{r}_1$ according
    to Proposition~\ref{proposition:HomogenesStarKap}, i.e.\
    $\Eulerfib \mathsf{r}_1 = \mathsf{r}_1$ together with the fact
    that $\degsst \mathsf{r}_1 = \mathsf{r}_1$ directly implies that
    $\deg_\lambda \mathsf{r}_1 = 0$, i.e.\ $\mathsf{r}_1$ is purely
    \emph{classical}. Furthermore, the above proposition shows that
    $r$ decomposes into one classical part $\mathsf{r}_1$ that is
    independent of $B$ and another non-classical part $\mathsf{r}_0$
    that completely encodes the dependence of $r$ on $B$.
\end{Remark}

To write down the explicit formulas for the $\starkap$-product with a
function $\pi^*u$ the following notions turn out to be very useful. By
$\SymD : \Gamma^\infty(\Sym^\bullet E^*) \longrightarrow \Gamma^\infty
(\Sym^{\bullet + 1}E^*)$ we denote the operator of \emph{symmetric
  $E$-covariant derivation}
\begin{equation*}
    \SymD = e^\alpha \nabla_{e_\alpha}.
\end{equation*}
For functions $u\in C^\infty(M)$ and one-forms $\alpha \in
\Gamma^\infty(E^*)$ we def\/ine the algebraic dif\/ferential operators
$\mathcal{F} (u)$ and $\mathcal{F} (\alpha)$ on $\S{\bullet}(E)$ by
\begin{equation*}
    \mathcal{F} (u) t = u t \qquad\textrm{and}\qquad \mathcal{F}( \alpha)
    t = \inss(\alpha)t
\end{equation*}
and extend the map $\mathcal{F}$ to a homomorphism from
$\S{\bullet}(E^*)=\bigoplus_{k=0}^\infty\Gamma^\infty(\Sym^k E^*)$ to
the algebraic dif\/ferential operators on $\S{\bullet}(E)$. For obvious
reasons we sometimes refer to the operators in the image of
$\mathcal{F}$ as \emph{fibre derivatives}, see
\cite{bordemann.neumaier.pflaum.waldmann:2003a} for a motivation.
Analogously, we def\/ine $\mathsf{F}(u)f = \pi^* u f$ and
$\mathsf{F}(\alpha) f = \Lie_{\alpha^\ver} f$ for $f \in
C^\infty(E^*)$ and extend $\mathsf{F}$ to a homomorphism from
$\S{\bullet}(E^*)$ to the dif\/ferential operators on $C^\infty(E^*)$.
Clearly, for all $\alpha \in \S{\bullet}(E^*)$ we therefore have the
relation $\mathcal{J}^{-1}\circ \mathsf{F}(\alpha) \circ \mathcal{J} =
\mathcal{F}(\alpha)$.
\begin{Proposition}
    \label{proposition:Sternkommutator}
    For all choices of $\D_E$-closed $B \in \Gamma^\infty(\Anti^2
    E^*)[[\lambda]]$ one has:
    \begin{enumerate}
    \itemsep=0pt
    \item\label{item:standardantistandard} The star products
        $\starstd$ resp.\ $\starastd$ are of \emph{standard-ordered
          type} resp.\ \emph{anti-standard-ordered} type in the sense
        that
        \begin{equation}
            \label{eq:standardantistandard}
            \pi^*u \starstd  f  = \pi^*u f\qquad\textrm{and}
            \qquad f \starastd \pi^*u = f \pi^*u
        \end{equation}
        for all $u \in C^\infty(M)[[\lambda]]$ and all $f \in
        C^\infty(E^*)[[\lambda]]$. Equivalently, $\starstd$ $($resp.\
        $\starastd)$ differentia\-tes the first $($resp.\ second$)$ argument
        only in fibre directions.
    \item\label{item:FedosovTaylorvonu} The Fedosov--Taylor series of
        $u \in C^\infty (M)[[\lambda]]$ is explicitly given by
        \begin{equation*}
            \tau_\kappa (u) = \exp (\SymD) u.
        \end{equation*}
    \item\label{item:linksrechtsMultiplikation} In addition, the
        $\starkap$-left and $\starkap$-right multiplications with
        $\pi^*u$ can be expressed as
        \begin{equation}
            \label{eq:linksrechtsMultiplikation}
            \pi^*u \starkap f  =
            \mathsf{F} \left(
                \exp\left(
                    \kappa \I \lambda \SymD
                \right) u
            \right) f
            \qquad\textrm{and}\qquad
            f \starkap \pi^*u =
            \mathsf{F} \left(
                \exp\left(
                    - (1-\kappa) \I\lambda \SymD
                \right) u
            \right) f.
        \end{equation}
        Consequently, for all $u \in C^\infty(M)[[\lambda]]$ one has
        \begin{equation*}
            \ad_{\starkap} (\pi^*u) =
            \mathsf{F}
            \left(
                \frac{\exp\left(\kappa \I\lambda\SymD\right)
                  - \exp\left(-(1-\kappa)\I\lambda\SymD\right)}
                {\SymD} \D_E u
            \right).
        \end{equation*}
    \item\label{item:glattaufMundeformed} For all $\starkap$ the
        functions $\pi^*C^\infty(M)[[\lambda]]$ are a subalgebra with
        $\starkap$ being the undeformed commutative product.
    \end{enumerate}
\end{Proposition}
\begin{proof}
    Assertion~\refitem{item:standardantistandard} follows from 
    $u \aststd t = \sigma (\taustd (u)\fibstd \taustd (t)) =
    \sigma ( (\sigma (\taustd(u))) \fibstd \taustd (t)) = \sigma (u
    \fibstd \taustd(t)) = \sigma (u \taustd(t))= u t$, where the
    second equality is due to the explicit shape of $\fibstd$. Observe
    that one does not need to know $\taustd(u)$ in order to prove this
    statement. Analogously one shows $t \astastd u = t u$ for all $u
    \in C^\infty(M)[[\lambda]]$ and all $t \in \S{}$. From the very
    def\/initions of $\starstd$ and $\starastd$ we can therefore
    conclude that the equations \eqref{eq:standardantistandard} hold
    since they hold on $\Pol^\bullet(E^*)[[\lambda]]$. The second
    statement follows using associativity.  For the proof of
    \refitem{item:FedosovTaylorvonu} we just have to show that $\exp
    (\SymD) u = u + \delta^{-1}\left( D \exp (\SymD) u +
        \frac{\I}{\lambda} \ad_\kappa(r)\exp (\SymD) u\right)$, which
    is straightforward observing that for all $a \in \W{}$ with $\degs
    a = k a$ we have $\delta^{-1} D a = \frac{1}{k +1} \SymD a$ and
    $\delta^{-1} \frac{\I}{\lambda}\ad_\kappa (r) a = \delta^{-1}
    (\inss (e^\alpha)r \inss(e_\alpha) a) = 0$, where the second to
    last equality is a consequence of the explicit shape of $\fibkap$
    and
    Proposition~\ref{proposition:rproperties}~\refitem{item:rsimpler}.
    Moreover, we have used that $\delta^{-1}r = 0$ according to
    \eqref{eq:requations}. The f\/irst statement in
    \refitem{item:linksrechtsMultiplikation} follows using
    \refitem{item:FedosovTaylorvonu} similarly to the computation in
    \refitem{item:standardantistandard}
    \begin{gather*}
        u \astkap t
        = \sigma (\tau_\kappa(u) \fibkap \tau_\kappa (t))
        \\
        \phantom{u \astkap t}{}
         = \sum_{k =0}^\infty \frac{1}{k!}
        \left(\kappa\I\lambda\right)^k
        \inss(e_{\alpha_1}) \cdots \inss (e_{\alpha_k})
        \frac{1}{k!}\SymD^k u
        \inss(e^{\alpha_1})\cdots \inss(e^{\alpha_k}) t \\
        \phantom{u \astkap t}{} = \mathcal{F} \left(
            \exp\left(\kappa\I\lambda\SymD\right)
            u
        \right) t.
    \end{gather*}
    From this equation and the relation between $\mathsf{F}$ and
    $\mathcal{F}$ we conclude the f\/irst equation in
    \eqref{eq:linksrechtsMultiplikation}. The proof of the second
    equation for the right-multiplication with $\pi^*u$ is completely
    analogous. The last assertion in
    \refitem{item:linksrechtsMultiplikation} then follows from the
    combination of the equations for $\pi^*u \starkap f$ and $f
    \starkap \pi^*u$ observing $\SymD u = \D_E u$.  Finally,
    \refitem{item:glattaufMundeformed} is a direct consequence of
    \refitem{item:linksrechtsMultiplikation}.
\end{proof}

We conclude this section noting that~-- in case $B$ has been chosen
appropriately~-- the star products $\starweyl$ are in fact of Weyl
type, i.e.\ the complex conjugation $\mathrm{C}$ as well as the
$\lambda$-parity operator $\mathrm{P} = (-1)^{\deg_\lambda}$ are
anti-automorphisms of $\starweyl$ (cf.\ \cite{neumaier:2002a}). Note
that we consider the formal parameter $\lambda$ as real and hence
def\/ine $\mathrm{C}\lambda = \lambda$.
\begin{Proposition}
    For $\starweyl = \star_{\frac{1}{2}}$ we have the following
    statements:
    \begin{enumerate}
    \itemsep=0pt
    \item\label{item:Hermitesch} In case $\mathrm{C} B = B$ the star
        product $\starweyl$ is \emph{Hermitian}, i.e.\ $\mathrm{C}( f
        \starweyl g ) = (\mathrm{C}g)\starweyl (\mathrm{C} f)$ for all
        $f, g\in C^\infty(E^*)[[\lambda]]$.
    \item\label{item:Paritaet} In case $\mathrm{P} B = B$ the star
        product $\starweyl$ has the \emph{$\lambda$-parity property},
        i.e.\ $\mathrm{P}( f \starweyl g ) = (\mathrm{P}g)\starweyl
        (\mathrm{P} f)$ for all $f,g\in C^\infty(E^*)[[\lambda]]$.
    \item\label{item:WeylTyp} In case $\mathrm{C} B = B = \mathrm{P}
        B$ the star product $\starweyl$ is of \emph{Weyl type}.
    \end{enumerate}
\end{Proposition}
\begin{proof}
    For the proof of \refitem{item:Hermitesch} we f\/irst note that
    $\mathrm{C}(a \fibweyl b) = (-1)^{k l} (\mathrm{C} b) \fibweyl
    (\mathrm{C} a)$ for all $a \in \WSL{}{}{k}$, $b \in \WSL{}{}{l}$.
    With $\mathrm{C}B = B$ and this property it is straightforward to
    show that $\mathrm{C}r$ satisf\/ies the same equations as $r$ and
    hence $\mathrm{C}r = r$ by the uniqueness of the solution of
    \eqref{eq:requations}. Therefore $\mathrm{C}$ commutes with
    $\FedDer_{\scriptscriptstyle\mathrm{Weyl}}$ implying
    $\mathrm{C}\tau_{\scriptscriptstyle\mathrm{Weyl}} (s) =
    \tau_{\scriptscriptstyle\mathrm{Weyl}} (\mathrm{C} s)$ for all $s
    \in \S{}$. From this equation and the behaviour of~$\fibweyl$
    with respect to complex conjugation the identity $\mathrm{C}(s \astweyl t)
    = (\mathrm{C} t) \astweyl (\mathrm{C} s)$ for all $s,t \in \S{}$
    is obvious implying the assertion. The proof
    of~\refitem{item:Paritaet} is completely analogous to that
    of~\refitem{item:Hermitesch} replacing $\mathrm{C}$ by~$\mathrm{P}$ and~\refitem{item:WeylTyp} is just the combination
    of~\refitem{item:Hermitesch} and~\refitem{item:Paritaet}.
\end{proof}

%
%

\section{Equivalence transformations}
\label{sec:EquivalenceTransformations}

In this section we will explicitly construct isomorphisms or
equivalence transformations, respectively, relating the star products
obtained for dif\/ferent values $\kappa$ and $\kappa'$ of the ordering
parameter, dif\/ferent torsion-free $E$-connections $\nabla$ and
$\nabla'$, and from cohomologous $E$-two-forms $B$ and $B'$.

We begin with the construction of isomorphisms between products
$\astkap$ and $\astStrkap$ obtained using the same ordering
parameter $\kappa$, the same $E$-connection $\nabla$ and dif\/ferent but
cohomologous formal series $B$ and $B'$ of $\D_E$-closed
$E$-two-forms. We look for an automorphism $\mathcal{A}_h$ of
$(\WSL{}{}{},\fibkap)$ of the form
\begin{equation}
    \label{eq:faserweiserAuto}
    \mathcal{A}_h = \exp\left(
        - \frac{\I}{\lambda} \adkap (h)\right) \qquad \textrm{with}\quad
    h \in \W{1} \quad \textrm{and} \quad \sigma (h)=0
\end{equation}
such that
\begin{equation*}
    \FedDer'_\kappa = \mathcal{A}_h\FedDer_\kappa \mathcal{A}_{-h}.
\end{equation*}
In case we can f\/ind such a f\/ibrewise automorphism the map
$\mathcal{I}_h : \S{} \longrightarrow \S{}$ which is def\/ined by
\begin{equation}
    \label{eq:BcohomologAequivalenzDef}
    \mathcal{I}_h (s) = \sigma\left(\mathcal{A}_h \tau_\kappa
        (s)\right)
\end{equation}
turns out to be an isomorphism from $\astkap$ to $\astStrkap$ since
the relation between the Fedosov derivations implies that
$\mathcal{A}_h \tau_\kappa(s) = \tau'_\kappa(\mathcal{I}_h s)$ holds
for all $s\in \S{}$. Now by a direct computation one f\/inds that
$\mathcal{A}_h\FedDer_\kappa \mathcal{A}_{-h} = \FedDer_\kappa +
\frac{\I}{\lambda} \adkap\left(\frac{\exp\left( - \frac{\I}{\lambda}
          \adkap(h) \right) - \id}{- \frac{\I}{\lambda} \adkap(h)}
    \FedDer_\kappa h \right)$ which is equal to $\FedDer'_\kappa$ if\/f
\begin{equation*}
    \frac{\I}{\lambda} \adkap\left(
        r - r' + \frac{\exp\left( - \frac{\I}{\lambda}
              \adkap(h) \right) - \id}{- \frac{\I}{\lambda} \adkap(h)}
        \FedDer_\kappa h
    \right)=0.
\end{equation*}
By Lemma~\ref{lemma:SuperCenter}, this is the case if\/f there is a
formal series of $E$-one-forms $A$ such that
\begin{equation*}
    r - r' + \frac{\exp\left( - \frac{\I}{\lambda}
          \adkap(h) \right) - \id}{- \frac{\I}{\lambda} \adkap(h)}
    \FedDer_\kappa h = - A.
\end{equation*}
Now using that $h$ is of $\Sym^\bullet E$-degree $0$ together with
Proposition~\ref{proposition:rproperties}~\refitem{item:rsimpler} it
is easy to show that $\FedDer_\kappa h$ is also of $\Sym^\bullet
E$-degree $0$ which implies that the above equation reduces to
\begin{equation}
    \FedDer_\kappa h = r' - r - A,
\end{equation}
since $\WSL{}{0}{}$ is a super-commutative subalgebra of
$(\WSL{}{}{},\fibkap)$. Applying $\FedDer_\kappa$ to this equation, we
get using the equations for $r$ and $r'$ that the necessary condition
for the solvability of this equation is that $A$ satisf\/ies
\begin{equation*}
    B - B' = \D_E A.
\end{equation*}
But this is also suf\/f\/icient since the $\FedDer_\kappa$-cohomology is
trivial on elements with positive $\Anti^\bullet E^*$-degree and the
solution $h_A$ is explicitly given by $h_A = \FedDer_\kappa^{-1}(r' -
r - A )$. We thus have shown the f\/irst part of the following:
\begin{Theorem}
    \label{theorem:Bcohomologequivalent}
    Let $\astkap$ and $\astStrkap$ be Fedosov products on $\S{}$
    obtained using the same $E$-connection and different but
    cohomologous $\D_E$-closed $E$-two-forms $B$ and $B'$.
    \begin{enumerate}
    \itemsep=0pt
    \item In case $A$ satisfies $B- B' = \D_E A$ there is a uniquely
        determined element $h_A \in \W{1}$ with $\sigma(h_A)=0$ such
        that $\FedDer_\kappa h_A = r' - r - A$ namely
        \begin{equation*}
            h_A = \FedDer_\kappa^{-1}(r' - r - A ).
        \end{equation*}
        With the so-constructed $h_A$ one has $\FedDer'_\kappa =
        \mathcal{A}_{h_A} \FedDer_\kappa \mathcal{A}_{- h_A}$ and thus
        $\mathcal{I}_A = \mathcal{I}_{h_A}$ according to equation
        \eqref{eq:BcohomologAequivalenzDef} defines an algebra
        isomorphism from $(\S{},\astkap)$ to $(\S{},\astStrkap)$.
    \item The above element $h_A \in \W{1}$ is explicitly given by
        \begin{equation*}
            h_A =   \frac{\exp(\SymD) - \id}{\SymD} A,
        \end{equation*}
        where the formal series of $E$-one-forms $A$ is considered as
        a one-form in the symmetric part of $\WSL{}{}{}$, i.e.\ as
        element of $\mathcal{W}^1$. Moreover, the isomorphism
        $\mathcal{I}_{A}$ takes the concrete form
        \begin{equation*}
            \mathcal{I}_{A} s = \exp \left(
                \mathcal{F}\left(
                    \frac{\exp\left(\kappa \I\lambda \SymD\right)
                      - \exp\left(-(1-\kappa) \I\lambda\SymD\right)}
                    {\I\lambda \SymD} A
                \right)
            \right)s
        \end{equation*}
        for $s \in \S{}$.
    \item The map $\mathsf{I}_{A}$ defined by
        \begin{equation}
            \label{eq:BcohomologIsoDef}
            \mathsf{I}_{A} =  \Phi_{A_0}^* \circ \exp \left(
                \mathsf{F}\left(
                    \frac{\exp\left(\kappa \I\lambda\SymD\right)
                      - \exp\left(-(1-\kappa) \I\lambda\SymD\right)}
                    {\I\lambda \SymD} A - A_0
                \right)
            \right),
        \end{equation}
        where $A_0$ denotes the classical part of $A$ and $\Phi_{A_0}:
        E^* \longrightarrow E^*$ denotes the fibre translating
        diffeomorphism defined by
        \begin{equation*}
            \Phi_{A_0}(\alpha_q) = \alpha_q + A_0(q),
        \end{equation*}
        provides an isomorphism from $(C^\infty(E^*)[[\lambda]],
        \starkap)$ to $(C^\infty(E^*)[[\lambda]],\starStrkap)$.
    \item In case $\deg_\lambda B = B$ and $\deg_\lambda B' = B'$, we
        can choose $A$ such that $\deg_\lambda A = A$, too. Then
        $\mathsf{I}_A$ is homogeneous, i.e.\ $\mathsf{I}_A \mathsf{H}
        = \mathsf{H} \mathsf{I}_A$.
    \end{enumerate}
\end{Theorem}
\begin{proof}
    For the proof of the second part we have to evaluate
    $\FedDer_\kappa^{-1}(r'-r-A)$ using the concrete shape of the
    homotopy operator $\FedDer_\kappa^{-1}$ explicitly. Since
    $[\delta^{-1},[\delta^{-1},D + \frac{\I}{\lambda} \adkap(r)]]=0$
    the terms involving $r' - r$ vanish due to
    $\delta^{-1}r=0=\delta^{-1}r'$ yielding $h_A = \frac{1}{\id -
      [\delta^{-1},D + \frac{\I}{\lambda} \adkap(r)]}A$, where $A$ is
    viewed as element in $\mathcal{W}^1$. But as $[\delta^{-1},D +
    \frac{\I}{\lambda} \adkap(r)]$ preserves $\W{}$ and
    $[\delta^{-1},\frac{\I}{\lambda} \adkap(r)]$ even vanishes on
    $\W{}$ we get $h_A = \frac{1}{\id - [\delta^{-1},D]}A =
    \frac{\exp(\SymD )- \id}{\SymD} A$ where in the last step we have
    used that $[\delta^{-1},D]a = \frac{1}{k+1}\SymD a$ for all $a \in
    \W{}$ with $\degs a = k a$ and have explicitly computed the
    geometric series. For the detailed computation of
    $\mathcal{I}_{A}s$ we f\/irst note that $\sigma(\adkap(a)b)=
    \sigma(\adkap(a) \sigma(b))$ for all $a \in \WSL{}{0}{}$ and all
    $b\in \WSL{}{}{}$ because of the explicit shape of $\fibkap$.
    Repeated use of this identity yields $\mathcal{I}_{A}s =
    \exp\left(- \frac{\I}{\lambda} \sigma \circ
        \adkap\left(\frac{\exp(\SymD) - \id}{\SymD}A \right) \right)
    s$. A further straightforward computation shows that
    \[
    \sigma\left(
        \adkap\left(
            \frac{\exp(\SymD) - \id}{\SymD}A
        \right)s
    \right)
    =
    \mathcal{F}\left(
        \frac{\exp\left(\kappa\I\lambda\SymD\right)
          - \exp\left(- (1-\kappa)\I\lambda\SymD\right)}
        {\SymD} A
    \right) s
    \]
    for all $s \in \S{}$ proving the formula for $\mathcal{I}_{A}$.
    Evidently, from the very def\/inition of $\starkap$ and
    $\starStrkap$ the map $\mathsf{I}_{A} = \mathcal{J} \circ
    \mathcal{I}_{A} \circ \mathcal{J}^{-1}$ def\/ines an isomorphism
    from $(\Pol^\bullet(E^*)[[\lambda]],\starkap)$ to
    $(\Pol^\bullet(E^*)[[\lambda]],\starStrkap)$. But as $A$ starts in
    order $0$ of the formal parameter the so-obtained map cannot be
    continued to a~well-def\/ined map on all of
    $C^\infty(E^*)[[\lambda]]$ without being suitably rewritten. This
    trouble is caused by the term $\exp(\mathcal{F}(A_0))$ which is
    perfectly well def\/ined on $\Pol^\bullet(E^*)$ but not on
    $C^\infty(E^*)$. But on $\Pol^\bullet(E^*)$ the map
    $\exp(\mathcal{F}(A_0))$ coincides with the pull-back
    $\Phi^*_{A_0}$ with the above f\/ibre translating dif\/feomorphism.
    Therefore on $\Pol^\bullet(E^*)[[\lambda]]$ the expression in
    equation \eqref{eq:BcohomologIsoDef} coincides with $\mathcal{J}
    \circ \mathcal{I}_{A}\circ \mathcal{J}^{-1}$ but it is moreover
    well def\/ined on all of $C^\infty(E^*)[[\lambda]]$. With the usual
    argument, that it suf\/f\/ices to verify identities relating the star
    products $\starkap$ and $\starStrkap$ on polynomial functions in
    order to prove them for all smooth functions, this implies that
    $\mathsf{I}_{A}$ def\/ines an isomorphism from
    $(C^\infty(E^*)[[\lambda]],\starkap)$ to
    $(C^\infty(E^*)[[\lambda]],\starStrkap)$. The last part is
    obvious.
\end{proof}

As a by-product of the above considerations and the formulas proven in
Proposition~\ref{proposition:Sternkommutator} we obtain a unique
characterization of the (quasi-inner) derivations and (inner)
self-equivalences of~$\starkap$ that are in the image of $\mathsf{F}$.
\begin{Corollary}
    The map
    \begin{equation*}
         A \mapsto
        \mathsf{F}\left(
            \frac{\exp\left(\kappa \I\lambda\SymD\right)
              - \exp\left(- (1-\kappa) \I\lambda \SymD\right)}
            {\I\lambda \SymD} A
        \right)
    \end{equation*}
    establishes a bijection between $\{ A \in
    \Gamma^\infty(E^*)[[\lambda]] \,|\, \D_E A =0\}$ and
    $\Der(C^\infty(E^*)[[\lambda]],\starkap) \cap \image(\mathsf{F})$.
    Since all self-equivalences of $\starkap$ are of the form
    $\exp\left(\I\lambda\mathrm{D}\right)$ with $\mathrm{D}
    \in\Der(C^\infty(E^*)[[\lambda]],\starkap)$ this also yields that
    all self-equivalences of $\starkap$ that lie in
    $\image(\mathsf{F})$ have the form
    \begin{equation*}
        \exp \left(
            \mathsf{F}\left(
                \frac{\exp\left(\kappa \I\lambda \SymD\right)
                  - \exp\left(-(1-\kappa) \I\lambda \SymD\right)}
                {\SymD} A
            \right)
        \right).
    \end{equation*}
    Moreover, the derivations obtained from $A = \D_E u$ with $u \in
    C^\infty(M)[[\lambda]]$ are of the form 
    $-\frac{\I}{\lambda}\ad_{\starkap}(\pi^* u)$ and hence are
    quasi-inner. Analogously, the self-equivalences corresponding to
    $A = \D_E u$ are given by $\exp(\ad_{\starkap}(\pi^* u))f =
    \E^{\pi^* u} \starkap f \starkap \E^{- \pi^*u}$ and hence are
    inner automorphisms.
\end{Corollary}

Let us now construct an equivalence transformation from $\starkap$ to
$\starStrkap$ obtained from dif\/ferent torsion-free $E$-connections
$\nabla$ and $\nabla'$ but coinciding formal series of $\D_E$-closed
$E$-forms $B$. We begin with the comparison of the corresponding maps
$D$ and $D'$ on $\WSL{}{}{}$. Evidently, $S_{\nabla - \nabla'}(s,t)=
\nabla_s t - \nabla'_s t$ def\/ines an element of $\Gamma^\infty(\Sym^2
E^* \otimes E)$ which we naturally can consider as $S_{\nabla -
  \nabla'} \in \mathcal{W}^2\otimes \mathcal{S}^1$. For $T_{\nabla -
  \nabla'} = \delta S_{\nabla - \nabla'}\in \mathcal{W}^1\otimes
\mathcal{S}^1 \otimes \Anti^1$ we f\/ind:
\begin{Lemma}
    With the definitions from above the following identities hold:
    \begin{enumerate}
    \itemsep=0pt
    \item\label{item:DDStrichDifferenz}
        \begin{equation*}
            D - D' = - \frac{\I}{\lambda} \adkap(T_{\nabla -\nabla'}).
        \end{equation*}
    \item\label{item:TNabNabStrProperties}
        \begin{equation*}
            \delta T_{\nabla - \nabla'}=0
            \qquad
            \textrm{and}
            \qquad
            D T_{\nabla - \nabla'}
            = R - R' - \frac{\I}{\lambda} T_{\nabla - \nabla'}
            \fibkap T_{\nabla - \nabla'}.
        \end{equation*}
\end{enumerate}
\end{Lemma}
\begin{proof}
    Part~\refitem{item:DDStrichDifferenz} follows from an easy
    computation using the def\/inition of $T_{\nabla - \nabla'}$ and
    $\fibkap$. The f\/irst statement in
    \refitem{item:TNabNabStrProperties} is trivial since $\delta^2=0$
    and the f\/irst formula involving $R$ and $R'$ follows from squaring
    the identity $D' = D + \frac{\I}{\lambda} \adkap(T_{\nabla -
      \nabla'})$ and using the def\/initions of $R$ and $R'$.
\end{proof}

To construct an equivalence from $\astkap$ to $\astStrkap$ we again
try to f\/ind a f\/ibrewise automorphism $\mathcal{A}_h$ as in
\eqref{eq:faserweiserAuto} such that $\FedDer'_\kappa =\mathcal{A}_h
\FedDer_\kappa \mathcal{A}_{-h}$ with $\sigma(h)=0$. As above this
equation turns out to be equivalent to
\begin{equation*}
    T_{\nabla - \nabla'} + r' - r - \frac{\exp\left( - \frac{\I}{\lambda}
          \adkap(h) \right) - \id}{- \frac{\I}{\lambda} \adkap(h)}
    \FedDer_\kappa h
\end{equation*}
being a central element in $(\WSL{}{}{},\fibkap)$. The following
proposition states the existence of such an element and hence yields
the desired equivalence from $\starkap$ to $\starStrkap$.
\begin{Proposition}\qquad{}
    \begin{enumerate}
    \itemsep=0pt
    \item\label{item:hconnectionExists} There is a uniquely determined
        element $h_{\nabla - \nabla'} \in \WS{2}{}$ with
        $\sigma(h_{\nabla - \nabla'})=0$ such that
        \begin{equation}
            \label{eq:hconnectionGlch}
            T_{\nabla - \nabla'} + r' - r - \frac{\exp\left( -
                  \frac{\I}{\lambda} \adkap(h_{\nabla-\nabla'}) \right) - \id}{-
              \frac{\I}{\lambda} \adkap(h_{\nabla-\nabla'})}
            \FedDer_\kappa h_{\nabla - \nabla'} = 0,
         \end{equation}
         which can be determined recursively from
         \begin{gather*}
             h_{\nabla - \nabla'} = \delta^{-1}\Bigg(
                 D h_{\nabla - \nabla'} + \frac{\I}{\lambda} \adkap(r)
                 h_{\nabla - \nabla'}\\
                 \phantom{h_{\nabla - \nabla'} =}{}
                 -  \frac{-
                   \frac{\I}{\lambda} \adkap(h_{\nabla- \nabla'})}{\exp\left( -
                       \frac{\I}{\lambda} \adkap(h_{\nabla - \nabla'})
                   \right) - \id} (T_{\nabla - \nabla'} + r' - r )
             \Bigg).
         \end{gather*}
         With the so-determined $h_{\nabla-\nabla'}$ the map
         $\mathcal{E}_{\nabla- \nabla'}$ def\/ined by
         \begin{equation*}
             \mathcal{E}_{\nabla- \nabla'} s = \sigma
             (\mathcal{A}_{h_{\nabla- \nabla'}} \tau_\kappa(s))
         \end{equation*}
         is an equivalence transformation from
         $(\S{}, \astkap)$ to $(\S{}, \astStrkap)$.
     \item\label{item:hconnectionsimpler} Moreover, $h_{\nabla -
           \nabla'}$ is of $\Sym^\bullet E$-degree at most $1$ and
         hence it satisfies the simpler recursion formula
        \begin{gather}
                        h_{\nabla - \nabla'} = \delta^{-1} \Bigg(
                D h_{\nabla - \nabla'} -
                \{ r, h_{\nabla - \nabla'}\}_{\fib}\nonumber\\
                \phantom{h_{\nabla - \nabla'} =}{}
                -  \frac{\{ h_{\nabla - \nabla'},
                  \cdot \}_{\fib}}{\exp\left( \{
                      h_{\nabla - \nabla'},\cdot \}_{\fib} \right) - \id}
                (T_{\nabla - \nabla'} + r' - r )
            \Bigg).\label{eq:hNetteRekursion}
         \end{gather}
     \item\label{item:hconnectionEquivalenz} $\mathcal{E}_{\nabla-
           \nabla'}$ induces an equivalence transformation from
         $(C^\infty(E^*)[[\lambda]],\starkap)$ to
         $(C^\infty(E^*)[[\lambda]],\starStrkap)$ by the unique
         extension from $\Pol^\bullet(E^*)[[\lambda]]$ to
         $C^\infty(E^*)[[\lambda]]$ of the equivalence transformation
         $\mathsf{E}_{\nabla- \nabla'} = \mathcal{J}\circ
         \mathcal{E}_{\nabla- \nabla'} \circ \mathcal{J}^{-1} :
         (\Pol^\bullet (E^*)[[\lambda]],\starkap)\longrightarrow
         (\Pol^\bullet(E^*)[[\lambda]],\starStrkap)$.
     \item In case $\deg_\lambda B = B$ we have $\mathcal{H} h_{\nabla
           - \nabla'} = h_{\nabla - \nabla'}$ and hence
         $\mathsf{E}_{\nabla - \nabla'}$ is homogeneous.
     \end{enumerate}
\end{Proposition}
\begin{proof}
    Solving equation \eqref{eq:hconnectionGlch} for $\FedDer_\kappa h$
    and applying $\FedDer_\kappa$ one obtains by a straightforward but
    cumbersome computation that the necessary condition for the
    solvability of this equation is satisf\/ied. But this also turns out
    to be suf\/f\/icient since \eqref{eq:DeltaHomotopies} applied to $h$
    yields the recursion formula for $h_{\nabla - \nabla'}$ which has
    a unique solution by the usual f\/ixed point argument. For the
    so-constructed $h_{\nabla - \nabla'}$ it is evident that
    $\mathcal{E}_{\nabla-\nabla'}$ is an isomorphism as stated.  Using
    Proposition~\ref{proposition:TauDifferOrd} and the fact that
    $h_{\nabla - \nabla'}\in \WS{2}{}$ it is lenghty but not dif\/f\/icult
    to show that $\mathcal{E}_{\nabla-\nabla'}$ is in fact a formal
    series of dif\/ferential operators that starts with the identity.
    The proof of \refitem{item:hconnectionsimpler} is a
    straightforward induction on the total degree using
    Proposition~\ref{proposition:rproperties}~\refitem{item:rsimpler},
    the shape of $T_{\nabla - \nabla'}$, and the fact that the
    $\fibkap$-super-commutator of two elements that are of
    $\Sym^\bullet E$-degree at most~$1$ is of $\Sym^\bullet E$-degree
    at most~$1$, too, and equals $\I\lambda$ times the f\/ibrewise
    Poisson bracket $\{\cdot,\cdot\}_{\fib}$.
    Assertion~\refitem{item:hconnectionEquivalenz} is a direct
    consequence of the fact that $\mathcal{E}_{\nabla -\nabla'}$ is a
    formal series of dif\/ferential operators which is an equivalence
    between dif\/ferential associative products on $\S{}$. The last part
    follows from
    Proposition~\ref{proposition:HomogenesStarKap}~\refitem{item:homogenBedingunganB}
    together with \eqref{eq:hNetteRekursion} since $\mathcal{H}
    T_{\nabla - \nabla'} = T_{\nabla - \nabla'}$.
\end{proof}

Finally, we can now compare the star products $\starkap$ and
$\starkapStr$ obtained from dif\/ferent ordering parameters $\kappa$ and
$\kappa'$ but identical $\nabla$ and $B$. Actually, here the procedure
is a little more involved since we do not only need a f\/ibrewise
automorphism but also the f\/ibrewise equivalence transformation
$\mathcal{M}_{\kappa' - \kappa}$ from $\fibkap$ to $\fibkapStr$
def\/ined in \eqref{eq:faserweiserNeumaier}. We consider
$\Hat{\mathcal{D}}_{\kappa'} = \mathcal{M}_{\kappa' - \kappa}
\FedDer_\kappa \mathcal{M}_{\kappa - \kappa'}$ which is evidently a
super-derivation of $\fibkapStr$ with square $0$. Using
Proposition~\ref{proposition:rproperties}~\refitem{item:rsimpler} and the
commutation relations for $\Delta_{\fib}$ we f\/ind
\begin{equation*}
    \Hat{\FedDer}_{\kappa'} = \FedDer_{\kappa'} + (\kappa - \kappa')
    \adkapStr(\Delta_{\fib} r)
\end{equation*}
which in general is evidently dif\/ferent from $\FedDer_{\kappa'}$.
Therefore, we try to f\/ind $\mathcal{A}'_h = \exp\left(-
    \frac{\I}{\lambda} \adkapStr(h)\right)$ with $h\in \W{2}$ and
$\sigma(h)=0$ such that $\mathcal{A}'_h \Hat{\FedDer}_{\kappa'}
\mathcal{A}'_{-h} = \FedDer_{\kappa'}$ since then analogously to the
considerations before we can explicitly def\/ine an equivalence
transformation from $(\S{},\astkap)$ to $(\S{}, \astkapStr)$.
\begin{Proposition}\qquad
    \label{proposition:kappaEquivalence}
    \begin{enumerate}
    \itemsep=0pt
    \item\label{item:hkappakappaStrExists} For every $E$-one-form
        $\Gamma$ with $\D_E \Gamma = - \Delta_{\fib} R$ and all
        $\kappa,\kappa'\in \mathbb{R}$ there is a uniquely determined
        element $h_{\Gamma,\kappa'-\kappa} \in \W{2}$ with
        $\sigma(h_{\Gamma,\kappa'-\kappa})=0$ such that
        \begin{equation}
            \label{eq:hkappakappaStr}
            \FedDer_{\kappa'} h_{\Gamma,\kappa'-\kappa}
            = - \I \lambda (\kappa' - \kappa) (\Gamma +
            \Delta_{\fib} r)
        \end{equation}
        namely $h_{\Gamma,\kappa'-\kappa} = - \FedDer^{-1}_{\kappa'}
        \left(\I\lambda (\kappa' - \kappa) (\Gamma + \Delta_{\fib} r)
        \right)$.
    \item\label{item:hkappakappaStrEquivalenzI} With the element
        $h_{\Gamma,\kappa'-\kappa}$ constructed
        in \refitem{item:hkappakappaStrExists} one has
        $\FedDer_{\kappa'} = \mathcal{A}'_{h_{\Gamma,\kappa'-\kappa}}
        \mathcal{M}_{\kappa'-\kappa} \FedDer_\kappa
        \mathcal{M}_{\kappa - \kappa'}
        \mathcal{A}'_{h_{\Gamma,\kappa'-\kappa}}$ and hence
        \begin{equation}
            \label{eq:NeumaierOperatorDef}
            \mathcal{N}_{\Gamma,\kappa' - \kappa} s = \sigma\left(
            \mathcal{A}'_{h_{\Gamma,\kappa' - \kappa}}
            \mathcal{M}_{\kappa' -\kappa}\tau_\kappa(s)
            \right)
        \end{equation}
        defines an equivalence transformation from $(\S{},\astkap)$ to
        $(\S{},\astkapStr)$.
    \item\label{item:hkappakappaStrEquivalenzII} Finally, the map
        $\mathsf{N}_{\Gamma,\kappa' - \kappa} = \mathcal{J} \circ
        \mathcal{N}_{\Gamma,\kappa' - \kappa} \circ \mathcal{J}^{-1}$
        defines an equivalence transformation from
        $(\Pol^\bullet(E^*)[[\lambda]], \starkap)$ to
        $(\Pol^\bullet(E^*)[[\lambda]], \starkapStr)$ that uniquely
        extends to an equivalence transformation from $(C^\infty (E^*)
        [[\lambda]], \starkap)$ to $(C^\infty (E^*)[[\lambda]],
        \starkapStr)$.
    \item\label{item:HomogeneousAgain} In case $\deg_\lambda B = B$,
        the equivalence transformation $\mathsf{N}_{\Gamma, \kappa' -
          \kappa}$ is homogeneous.
    \end{enumerate}
\end{Proposition}
\begin{proof}
    Using the shape of $\mathcal{A}'_h$ which is a f\/ibrewise
    automorphism of $\fibkapStr$ we f\/ind that $\mathcal{A}'_h
    \Hat{\FedDer}_{\kappa'} \mathcal{A}'_{-h}$ $ = \FedDer_{\kappa'}$ if\/f
    \[
        \frac{\exp\left(- \frac{\I}{\lambda}\adkapStr(h) \right) -
          \id}{- \frac{\I}{\lambda}\adkapStr(h)} \FedDer_{\kappa'} h
        + \I \lambda (\kappa' - \kappa) \mathcal{A}'_h \Delta_{\fib}r
    \]
    is central. But now $\FedDer_{\kappa'} h \in \WSL{2}{0}{1}$ and
    $\Delta_{\fib} r \in \WSL{}{0}{1}$ together with the fact that
    $\WSL{}{0}{}$ is a super-commutative subalgebra of
    $(\WSL{}{}{},\fibkapStr)$ imply that this is equivalent~to
    \[
    \FedDer_{\kappa'} h
    + \I\lambda (\kappa' - \kappa)\Delta_{\fib} r = A,
    \]
    with a formal series of $E$-one-forms $A$. Clearly, the necessary
    condition for the solvability of the latter equation is
    $\FedDer_{\kappa'}\left(\I\lambda (\kappa' - \kappa)\Delta_{\fib}
        r \right) = \D_E A$. In order to analyze this equation one
    applies $\mathcal{M}_{\kappa' -\kappa}$ to the equation solved by
    $r$ with the f\/ibrewise product $\fibkap$ and f\/inds 
    $\FedDer_{\kappa'}\left(\I\lambda (\kappa' - \kappa)\Delta_{\fib}
        r \right) = \I\lambda (\kappa' -\kappa)\Delta_{\fib}R$. Note
    that here the fact that $r$ does \emph{not} depend on $\kappa$
    proven in Proposition~\ref{proposition:rproperties} enters
    crucially. But since there is an $E$-one-form $\Gamma$ with $\D_E
    \Gamma= - \Delta_{\fib} R$ we can satisfy the above condition
    choosing $A = - \I\lambda(\kappa' -\kappa) \Gamma$. Now again the
    homotopy operator $\FedDer_{\kappa'}^{-1}$ permits to write down
    the unique solution $h_{\Gamma,\kappa' -\kappa}$ of
    \eqref{eq:hkappakappaStr}. The assertion of
    \refitem{item:hkappakappaStrEquivalenzI} is then satisf\/ied by
    construction of $h_{\Gamma,\kappa' -\kappa}$ and
    part~\refitem{item:hkappakappaStrEquivalenzII} follows again by
    the argument that $\mathcal{N}_{\Gamma,\kappa' -\kappa}$ is a
    formal series of dif\/ferential operators starting with the
    identity. For the last part, note that $\Delta_{\fib} r =
    \Delta_{\fib} \mathsf{r}_1$ is \emph{classical} whence the
    right-hand side of \eqref{eq:hkappakappaStr} is homogeneous of
    degree one with respect to $\mathcal{H}$. This implies that
    $h_{\Gamma, \kappa'-\kappa}$ is homogeneous of degree one, too.
    Since $\mathcal{M}_{\kappa' - \kappa}$ commutes with
    $\mathcal{H}$, the result follows.
\end{proof}
\begin{Remark}
    \label{remark:InEquivalence}
    From the above statements it follows that the only possible source
    for non-equi\-valent star products in our framework are non-trivial
    cohomology classes $[B - B'] \in \mathrm{H}^2_E(M)[[\lambda]]$.
    Note however, that even for $[B] \ne [B']$ the resulting star
    products can be equivalent: let $E = \mathfrak{g}$ be an
    \emph{Abelian} even-dimensional Lie algebra, viewed as Lie
    algebroid over $\{\mathrm{pt}\}$. Choose $B_0$ to be a
    non-degenerate element in $\Anti^2 \mathfrak{g}^*$, whence the
    Poisson bracket $\{ \cdot, \cdot\}_{B_0}$ is symplectic. Then any
    $B = B_0 + \lambda B_1$ and $B' = B_0 + \lambda B_1'$ are
    $\D_E$-closed and cohomologous only if $B_1 = B_1'$. Choosing $B_1
    \ne B_1'$, the star products are still equivalent, as all
    symplectic star products on a vector space are equivalent.
\end{Remark}

%
%

\section{Relation to the universal enveloping algebra}
\label{sec:Representations}

Now we relate the deformed algebra of polynomial functions
$\left(\Pol^\bullet(E^*)[\lambda], \starkap\right)$ for $B = 0$,
viewed as $\mathbb{C}[\lambda]$-subalgebra thanks to
Proposition~\ref{proposition:HomogeneousStarProducts}, to the
universal enveloping algebra $\mathcal{U}(E)$ of the Lie algebroid
$E$.  First we recall the def\/inition of $\mathcal{U}(E)$ according to
Rinehart \cite[Section~2]{rinehart:1963a}: The
$\mathbb{C}[\lambda]$-module $\left(C^\infty(M) \oplus
    \Gamma^\infty(E)\right)[\lambda]$ is a Lie algebra over
$\mathbb{C}[\lambda]$ via
\begin{equation*}
    \left[(u, s), (v, t)\right]
    = - \I\lambda \left( \varrho(s)v - \varrho(t)u, [s, t]_E\right).
\end{equation*}
Note that we have incorporated artif\/icially the pre-factor $\I\lambda$
which will allow for easier comparison.  In its (Lie algebraic)
universal enveloping algebra $\mathfrak{U}\left(\left(C^\infty(M)
        \oplus \Gamma^\infty(E)\right)[\lambda]\right)$ over
$\mathbb{C}[\lambda]$ one considers the two-sided ideal $\mathfrak{I}$
generated by the relations $(u, 0) \diamond (0,s) - (0, us)$ and $(u,
0) \diamond (v,0) - (uv, 0)$, where $\diamond$ denotes the product of
the universal enveloping algebra and $u, v \in C^\infty(M)[\lambda]$
and $s \in \Gamma^\infty(E)[\lambda]$. Then one def\/ines
\begin{equation*}
    \mathcal{U}(E) =
    \mathfrak{U}\left(
        \left(C^\infty(M) \oplus \Gamma^\infty(E)\right)[\lambda]
    \right)
    \big/ \mathfrak{I}.
\end{equation*}
The universal enveloping algebra $\mathcal{U}(E)$ (now in the Lie
algebroid sense) is still f\/iltered by the number of factors from
$\Gamma^\infty(E)$. We denote by
$\mathfrak{Gr}^\bullet(\mathcal{U}(E))$ the corresponding graded
$\mathbb{C}[\lambda]$-module, i.e.\ $\mathfrak{Gr}^k(\mathcal{U}(E)) =
\mathcal{U}^k(E) \big/ \mathcal{U}^{k-1}(E)$, endowed with its
canonical algebra structure. It turns out that
$\mathfrak{Gr}^\bullet(\mathcal{U}(E))$ is a commutative unital
$\mathbb{C}[\lambda]$-algebra, which is immediate from the particular
form of the relations used in the def\/inition of $\mathcal{U}(E)$.
This algebra is still generated by $C^\infty(M)[\lambda]$ and
$\Gamma^\infty(E)[\lambda]$ and $C^\infty(M)[\lambda] \longrightarrow
\mathfrak{Gr}^0(\mathcal{U}(E))$ turns out to be an algebra
morphism. Thus $\mathfrak{Gr}^\bullet(\mathcal{U}(E))$ is a
$C^\infty(M)[\lambda]$-algebra, generated by
$\Gamma^\infty(E)[\lambda]$. Hence we have a canonical algebra
morphism
\begin{equation}
    \label{eq:PBW}
    \Sym^\bullet_{C^\infty(M)[\lambda]}
    \left(\Gamma^\infty(E)[\lambda]\right)
    =
    \mathcal{S}^\bullet(E)[\lambda]
    \longrightarrow
    \mathfrak{Gr}^\bullet(\mathcal{U}(E))
\end{equation}
by the universal property of the symmetric algebra. Since
$\Gamma^\infty(E)[\lambda]$ is a projective
$C^\infty(M)[\lambda]$-module, which follows directly from the
Serre--Swan theorem, it follows from \cite[Theorem~3.1]{rinehart:1963a}
that \eqref{eq:PBW} is in fact an isomorphism, i.e.\ one has a
Poincar\'e--Birkhof\/f--Witt like theorem. In particular, we can identify
$u \in C^\infty(M)[\lambda]$ and $s \in \Gamma^\infty(E)[\lambda]$
with their images in $\mathcal{U}(E)$.

Now consider $\mathcal{S}^\bullet(E)[\lambda]$ (or equivalently,
$\Pol^\bullet(E^*)[\lambda]$) with the (star) product $\astweyl$ (or
$\starweyl$), where we have to choose $B = 0$.
\begin{Proposition}
    \label{proposition:Covariance}
    For $B = 0$ we have
    \begin{equation}
        \label{eq:WeylCommutators}
        [u, v]_{\astweyl} = 0,
        \qquad
        [u, s]_{\astweyl} = \I \lambda \varrho(s)u,
        \qquad
        \textrm{and}
        \qquad
        [s, t]_{\astweyl} =  - \I \lambda [s, t]_E
    \end{equation}
    for $u, v \in C^\infty(M)[\lambda]$ and $s, t \in
    \Gamma^\infty(E)[\lambda]$.
\end{Proposition}
\begin{proof}
    The f\/irst two statements are true by homogeneity, the last follows
    from homogeneity and the Weyl type property.
\end{proof}
\begin{Remark}
    \label{remark:Covariance}
    This proposition, transferred to $\Pol^\bullet(E^*)$, can be viewed
    as a covariance property of $\starweyl$ under the Lie algebra
    action of $C^\infty(M) \rtimes \Gamma^\infty(E)$ acting on
    $C^\infty(E^*)$ via inner Poisson derivations. Thus the map
    $\mathcal{J}$, restricted to $C^\infty(M) \rtimes
    \Gamma^\infty(E)$, can be considered as a quantum momentum map.
\end{Remark}

\begin{Lemma}
    \label{lemma:MorphismToUE}
    There exists a unique surjective unital
    $\mathbb{C}[\lambda]$-linear algebra morphism
    \begin{equation*}
        \phi: \mathcal{U}(E) \longrightarrow
        \left(\mathcal{S}^\bullet(E)[\lambda], \astweyl\right)
    \end{equation*}
    with $\phi(u) = u$ and $\phi(s) = s$ for $u \in
    C^\infty(M)[\lambda]$ and $s \in \Gamma^\infty(E)[\lambda]$.
\end{Lemma}
\begin{proof}
    By \eqref{eq:WeylCommutators} and the universal property of
    $\mathcal{U}(E)$ existence and uniqueness of such a morphism are
    clear. The surjectivity follows since
    $\mathcal{S}^\bullet(E)[\lambda]$ is generated by
    $C^\infty(M)[\lambda]$ and $\Gamma^\infty(E)[\lambda]$ by
    Proposition~\ref{proposition:HomogeneousStarProducts}.
\end{proof}

Note also that $\mathcal{S}^\bullet(E)[\lambda]$ is no longer graded
with respect to $\astweyl$ but only f\/iltered.
\begin{Lemma}
    \label{lemma:GrWeylstarIsSymmetric}
    The graded algebra
    $\mathfrak{Gr}^\bullet\left(\mathcal{S}^\bullet(E)[\lambda],
        \astweyl\right)$ is canonically isomorphic to the symmetric
    algebra $\mathcal{S}^\bullet(E)[\lambda]$.
\end{Lemma}
\begin{proof}
    This is an immediate consequence of
    Proposition~\ref{proposition:HomogeneousStarProducts}.
\end{proof}
\begin{Proposition}
    \label{proposition:Universal}
    $\left(\mathcal{S}^\bullet(E)[\lambda], \astweyl\right)$ is
    isomorphic to $\mathcal{U}(E)$ via $\phi$.
\end{Proposition}
\begin{proof}
    It remains to show that $\phi$ is injective. Here we can rely on
    the complete symbol calculus for Lie algebroids as established in
    \cite[Theorem~3]{nistor.weinstein.xu:1999a} which is (as our
    construction) based on the choice of a connection $\nabla$.
\end{proof}

{\samepage \begin{Remark}\label{remark:Universal}\qquad
        \begin{enumerate}
    \itemsep=0pt
    \item Since $\astweyl$ converges trivially on
        $\mathcal{S}^\bullet(E)[\lambda]$ for any value $\lambda =
        \hbar \in \mathbb{C}$, we can also pass to the universal
        enveloping algebra of the Lie algebroid where $\lambda$ is
        replaced by $\hbar$. This f\/inally establishes the contact to
        the work of Nistor, Weinstein, and Xu
        \cite{nistor.weinstein.xu:1999a}.
    \item For the case $E = \mathfrak{g}$, the above construction
        reproduces the results of Gutt \cite{gutt:1983a}.
    \item For $\kappa \ne \frac{1}{2}$, we also obtain isomorphisms of
        $(\mathcal{S}^\bullet(E)[\lambda], \astkap)$ with
        $\mathcal{U}(E)$ as all $\kappa$-ordered products are
        equivalent via the explicit equivalence transformation
        \eqref{eq:NeumaierOperatorDef}, which preserves the subspace
        $\mathcal{S}^\bullet(E)[\lambda] \subseteq
        \mathcal{S}^\bullet(E)[[\lambda]]$ thanks to
        Proposition~\ref{proposition:kappaEquivalence}~\refitem{item:HomogeneousAgain}.
    \end{enumerate}
\end{Remark}}

\begin{Remark}[PBW Approach]
    \label{remark:AlternativeConstruction}
    One may wonder whether a choice of a PBW isomorphism from the
    symmetric algebra $\mathcal{S}^\bullet(E)[\lambda]$ to the graded
    algebra $\mathfrak{Gr}^\bullet (\mathcal{U}(E))$ can be used to
    give an alternative construction of the star product $\astweyl$.
    This construction was proposed to us by one of the referees.  In
    fact, a posteriori Proposition~\ref{proposition:Universal} shows
    that this is possible by using the isomorphism $\phi$ and its
    inverse described above.  However, in this approach one would end
    up with a deformation quantization of the symmetric algebra
    $\mathcal{S}^\bullet(E)[\lambda]$, even homogeneous, of which it
    is a priori \emph{not} clear whether it will extend to a
    \emph{differential} star product on $C^\infty(E^*)[[\lambda]]$.
    The point is that the choice of $\phi$ can also be made
    dif\/ferently thereby immediately destroying the property that the
    resulting product is bidif\/ferential. Ultimately, the reason why it
    works with~$\phi$ is a rather involved combinatorial problem: in
    the Lie algebra case studied in \cite{gutt:1983a} a detailed
    analysis of the BCH series was necessary to prove this.
\end{Remark}
\begin{Example}[Non-dif\/ferential PBW isomorphism]
    In \cite{bordemann.neumaier.waldmann:1998a} an explicit counter
    example was given in the case of a cotangent bundle $E^* =
    T^*\mathbb{R}$ of the real line of a linear isomorphism between
    $\mathcal{S}(E^*)[\lambda]$ and $\mathcal{U}(E)$ compatible with
    the f\/iltering but not leading to a bidif\/ferential star product: in
    this case, the universal enveloping algebra is the algebra of
    dif\/ferential operators on~$\mathbb{R}$ and a candidate for a
    linear isomorphism like $\phi^{-1}$ is a \emph{quantization map}
    \[
    \varrho: \ \Pol(T^*\mathbb{R})
    \longrightarrow \mathrm{Dif\/fop}(\mathbb{R}),
    \]
    which is completely specif\/ied by its values on $\varrho(\chi p^k)$
    for $\chi \in C^\infty(\mathbb{R})$ and $k \in \mathbb{N}_0$ where
    $p$ denotes the coordinate along the f\/ibers. Then using
    $\varrho(\chi p^k)
    = \left(\frac{\lambda}{\I}\right)^k\chi
    \frac{\partial^k}{\partial x^k}$ for $k \ne 2$ but
    \[
    \varrho(\chi p^2)
    = \chi \left(
        - \lambda^2
        \frac{\partial^2}{\partial x^2} + \lambda
        \frac{\partial}{\partial x}
    \right)
    \]
    yields an isomorphism such that the pulled back product is
    \emph{not} bidif\/ferential (but still homogeneous). In view of this
    counter example the property of $\astweyl$ to be bidif\/ferential is
    far from being trivial. In fact, we are not aware of a simple
    proof of this fact, even for the above isomorphism~$\phi$.
\end{Example}

In view of this example it remains challenging to show that one can
obtain a bidif\/ferential star product for a particular (and
geometrically motivated) choice of the PBW isomorphism.  One way to
prove it independently of our construction is to use the
pseudo-dif\/ferential operator approach of Nistor, Weinstein and Xu: in
fact they discuss this construction in
\cite[Theorem~4]{nistor.weinstein.xu:1999a} on the level of polynomial
functions. However, they do not prove explicitly that the resulting
star product is actually bidif\/ferential. But using their integral
formulas of the pseudo-dif\/ferential calculus one should be able to
proceed as in the case of the tangent bundle and ordinary
pseudo-dif\/ferential operators. In any case, one should note that their
approach is much more analytic than our Fedosov construction: in fact
there are several possible generalizations of the Fedosov construction
to a completely algebraic framework not referring to underlying
geometry at all.
\begin{Remark}
    \label{remark:Extension}
    As it was mentioned to us by one of the referees, the
    $\D_E$-closed two-form $B$ can be used to obtain a central
    extension of the Lie algebroid $E$. For this Lie algebroid, say~$\tilde{E}$, one can use the construction for $B = 0$ and take an
    appropriate quotient later to get the deformation quantization of
    the polynomials $\mathcal{S}^\bullet(E)$. However, here it seems
    to be even more dif\/f\/icult to control whether the result is
    bidif\/ferential or not. In fact, geometrically the quotient would
    correspond to a restriction to a submanifold for which it is known
    that while Poisson structures might be tangential it can happen
    that no (bidif\/ferential) star product restricts, see the example
    in \cite{cahen.gutt.rawnsley:1996a}. Thus for this case, it seems
    to be unavoidable to use the above Fedosov machinery to obtain a
    bidif\/ferential star product.
\end{Remark}

Again, also in this case one might apply techniques from
\cite{nistor.weinstein.xu:1999a} in so far as the pseudo-dif\/ferential
calculus may be extended to sections of certain line bundles over the
base manifold. The resulting star product for the symbols should then
incorporate the Chern form of this line bundle as $B$. Of course, this
would limit the approach to integral two-forms. For the cotangent
bundle case this was shown to work in
\cite{bordemann.neumaier.pflaum.waldmann:2003a}. We plan to
investigate this relation to representation theory of the deformed
algebras in some following work including the corresponding Morita
theory.

In view of the above remarks one may wonder why one should impose the
condition that a star product is dif\/ferential at all. In fact, there
are interesting examples of ``star products'' def\/ined on polynomial
functions on a cotangent bundle which are \emph{not} dif\/ferential but
enjoy other nice feature, see e.g.\ the projectively equivariant
quantizations in \cite{duval.elgradechi.ovsienko:2004a}. However, many
constructions in deformation quantization rely heavily on the fact
that star products are local (and hence locally bidif\/ferential) like
e.g.\ the \v{C}ech cohomological approach to the classif\/ication
\cite{gutt.rawnsley:1999a, neumaier:2002a} and to the existence and
normalization of traces \cite{gutt.rawnsley:2002a} to mention just a
few.

\section{The trace in the unimodular case}
\label{sec:Trace}

We shall now construct a trace functional for any homogeneous star
product on $E^*$, in particular for $\starkap$ with the choice $B =
\lambda B_1$ with $B_1 \in \Gamma^\infty(\Anti^2 E^*)$, $\D_E B_1 =
0$, for the case of a \emph{unimodular} Lie algebroid.  To this end we
f\/irst have to recall some results from \cite{evens.lu.weinstein:1999a}
on the modular class, see also \cite{weinstein:1997a,
  kosmann-schwarzbach.weinstein:2005a,
  kosmann-schwarzbach.laurent-gengoux.weinstein:2007a} as well as
\cite{huebschmann:1999a} for a more algebraic approach.  However, we
shall use a slightly dif\/ferent presentation.

For a Lie algebroid $E$ there are essentially two ways of def\/ining a
divergence of $s \in \Gamma^\infty(E)$: if $\mu \in
\Gamma^\infty(|\Anti^n|T^*M)$ is a positive density on $M$ then we can
def\/ine
\begin{equation*}
    \divergence_\mu(s)
    = \frac{1}{\mu} \Lie_{\varrho(s)} \mu \in C^\infty(M).
\end{equation*}
If $\nu \in \Gamma^\infty(|\Anti^N|E)$ is a positive $E$-density, we
analogously can def\/ine
\begin{equation*}
    \divergence_\nu(s) = \frac{1}{\nu} \Lie^{\,E}_s \nu \in C^\infty(M),
\end{equation*}
where we use the fact that the $E$-Lie derivative also acts on the
$E$-densities in the usual way, see
e.g.\ \cite[Section~2.2.5]{waldmann:2006a:prebook} for an elementary
introduction to the calculus of densities. Alternatively, and probably
more familiar, one can use volume forms instead of densities, provided
$M$ as well as $E$ are orientable. For $u, v \in C^\infty(M)$ and $s,
t \in \Gamma^\infty(E)$ we have the usual identities
\begin{gather}
    \label{eq:divmuPropertiesI}
    \divergence_\mu(us) = u \divergence_\mu(s) + \varrho(s)u, \\
    \label{eq:divmuPropertiesII}
    \divergence_{\E^v\mu} (s) = \divergence_\mu(s) + \varrho(s) v, \\
    \label{eq:divmuPropertiesIII}
    \divergence_\mu([s, t]_E) = \varrho(s) \divergence_\mu(t) -
    \varrho(t) \divergence_\mu(s),
\\
    \label{eq:divnuPropertiesI}
    \divergence_\nu(us)  = u \divergence_\nu(s) - \varrho(s)u, \\
    \label{eq:divnuPropertiesII}
    \divergence_{\E^v\nu} (s)  = \divergence_\nu(s) + \varrho(s) v, \\
    \label{eq:divnuPropertiesIII}
    \divergence_\nu([s, t]_E)  = \varrho(s) \divergence_\nu(t) -
    \varrho(t) \divergence_\nu(s).
\end{gather}
Note the dif\/ferent signs in \eqref{eq:divmuPropertiesI} and
\eqref{eq:divnuPropertiesI}: this implies that the map
\begin{equation*}
    s \mapsto \tr\ad(s) = \divergence_\mu(s) + \divergence_\nu(s)
\end{equation*}
is $C^\infty(M)$-linear and hence $\tr\ad \in
\Gamma^\infty(E^*)$. From \eqref{eq:divmuPropertiesIII} and
\eqref{eq:divnuPropertiesIII} we see that $\D_E \tr\ad = 0$ and from
\eqref{eq:divmuPropertiesII} and \eqref{eq:divnuPropertiesII} it
follows that the dependence of $\tr\ad$ on $\mu$ and $\nu$ is only via
$\D_E$-exact terms whence $[\tr\ad] \in \mathrm{H}^1_E(M)$ is a
well-def\/ined class \cite{evens.lu.weinstein:1999a}:
\begin{Definition}[Modular class]
    \label{definition:ModularClass}
    The class $[\tr\ad] \in \mathrm{H}_E^1(M)$ is called the modular
    class of $E$ and $E$ is called unimodular if $[\tr\ad] = 0$.
\end{Definition}
In fact, $\tr\ad$ can be seen as the `trace of the adjoint
representation', even though in general the adjoint representation
`$\ad$' is not really def\/ined in the sense of a Lie algebroid
representation.  In the unimodular case we can f\/ind $\mu$ and $\nu$
such that $\tr\ad(s) = 0$ for all $s \in \Gamma^\infty(E)$.

In \cite{evens.lu.weinstein:1999a} the modular class of a Poisson
manifold $(M, \theta)$ was related to the modular class of the
corresponding Lie algebroid $T^*M$ induced by $\theta$ as in
Example~\ref{example:LieAlgebroids}~\refitem{item:PoissonLieAlgebroid}.
Here we shall give yet another interpretation of $\tr\ad$ in terms of
the Poisson manifold $E^*$. To this end we f\/irst construct a density
on $E^*$ out of $\mu$ and $\nu$ in the following way. With respect to
local coordinates $(q^1, \ldots, q^n, p_1, \ldots, p_N)$ induced by
coordinates $x^1, \ldots, x^n$ on $M$ and local basis sections $e_1,
\ldots, e_N$ of $E$, we def\/ine
\begin{equation*}
    (\mu\otimes\nu)\left(
        \frac{\partial}{\partial q^1}, \ldots,
        \frac{\partial}{\partial q^n},
        \frac{\partial}{\partial p_1}, \ldots,
        \frac{\partial}{\partial p_N}
    \right)
    =
    \pi^* \mu\left(
        \frac{\partial}{\partial x^1}, \ldots,
        \frac{\partial}{\partial x^n}
    \right)
    \pi^* \nu\left(
        e^1, \ldots, e^N
    \right).
\end{equation*}
The following lemma is a simple computation using the transformation
properties of densities.
\begin{Lemma}
    \label{lemma:muotimesnuDef}
    $\mu \otimes \nu \in \Gamma^\infty(|\Anti^{n+N}|T^*E^*)$ is a
    globally well-defined positive density on $E^*$.
\end{Lemma}
\begin{Lemma}
    \label{lemma:ConstantDensities}
    Let $\Omega \in \Gamma^\infty(|\Anti^{n+N}|T^*E^*)$ be a positive
    density on $E^*$. Then the following statements are equivalent:
    \begin{enumerate}
    \itemsep=0pt
    \item\label{item:OmegaConstant} $\Omega$ is constant along the
        fibres.
    \item\label{item:OmegaTensor} There exist positive densities $\mu
        \in \Gamma^\infty(|\Anti^n|T^*M)$ and $\nu \in
        \Gamma^\infty(|\Anti^N|E)$ such that $\Omega = \mu \otimes
        \nu$.
    \item\label{item:OmegaHomogeneous} $\Lie_\xi \Omega = N \Omega$.
    \end{enumerate}
\end{Lemma}
\begin{proof}
    Assume \refitem{item:OmegaConstant} and choose some positive
    densities $\tilde{\mu} \in \Gamma^\infty(|\Anti^n|T^*M)$ and $\nu
    \in \Gamma^\infty(|\Anti^N|E)$. Then we can write $\Omega = f
    \tilde{\mu} \otimes \nu$ with $f \in C^\infty(E^*)$ and $\Omega$
    is constant along the f\/ibres.  Thus $f = \pi^* u$ with $u \in
    C^\infty(M)$ and $\mu = u\tilde{\mu}$ and $\nu$ will fulf\/ill
    \refitem{item:OmegaTensor}. If \refitem{item:OmegaTensor} holds
    then \refitem{item:OmegaHomogeneous} is a simple computation. Now
    assume \refitem{item:OmegaHomogeneous} and choose some $\mu$ and
    $\nu$ as before. Since clearly $\Lie_\xi(\mu \otimes \nu) = N \mu
    \otimes \nu$ the \emph{function} $f = \Omega / (\mu \otimes \nu)$
    satisf\/ies $\Lie_\xi f = 0$. Hence $f$ and thus $\Omega$ are
    constant along the f\/ibres.
\end{proof}
The next proposition is a simple computation using e.g.\ the local
formulas for the Hamiltonian vector f\/ield
\eqref{eq:HamiltonianVectorField}:
\begin{Proposition}
    \label{proposition:LieXJsmuotimesnu}
    Let $\mu \in \Gamma^\infty(|\Anti^n|T^*M)$ and $\nu \in
    \Gamma^\infty(|\Anti^N|E)$ be positive densities. Then we have for
    $s \in \Gamma^\infty(E)$ and $u \in C^\infty(M)$
    \begin{gather}
        \label{eq:LieXJsmunu}
        \Lie_{X_{\mathcal{J}(s)}} (\mu \otimes \nu)
        = \pi^*(\tr\ad(s)) \mu \otimes \nu,
    \\
        \label{eq:LieXumunu}
        \Lie_{X_{\pi^*u}} (\mu \otimes \nu) = 0.
    \end{gather}
\end{Proposition}

We shall now interpret this proposition in terms of the modular
vector f\/ield on $E^*$.  Recall that the \emph{vertical lift}
\begin{equation*}
    \Gamma^\infty(\Anti^kE^*) \ni \omega
    \; \mapsto \;
    \omega^\ver \in \Gamma^\infty(\Anti^k TE^*)
\end{equation*}
induces a map in cohomology
\begin{equation}
    \label{eq:VerticalLiftCohomology}
    \mathrm{H}^\bullet_E(M) \ni [\omega]
    \; \mapsto \;
    [\omega^\ver] \in \mathrm{H}^\bullet_{\theta_E}(E^*),
\end{equation}
since we have (with our sign conventions) the relation
\begin{equation*}
    \left(\D_E \omega\right)^\ver = - \D_{\theta_E} \omega^\ver,
\end{equation*}
see e.g.\ \cite[Theorem~15]{grabowski.urbanski:1997a}.  Here
$\mathrm{H}^\bullet_{\theta_E}(E^*)$ denotes the Poisson cohomology
and $\D_{\theta_E} = \Schouten{\theta_E, \cdot}$.  In general, the map
\eqref{eq:VerticalLiftCohomology} is far from being surjective, the
trivial Lie algebroid provides a simple counter-example. However, it
turns out to be injective in general:
\begin{Theorem}
    \label{theorem:VerticalLiftInjective}
    The vertical lift induces an injective
    map \eqref{eq:VerticalLiftCohomology} in cohomology.
\end{Theorem}
\begin{proof}
    Let $\Phi_t(\alpha_q) = \E^t \alpha_q$ be the f\/low of the Euler
    vector f\/ield $\xi$ on $E^*$. For a $k$-vector f\/ield $X \in
    \Gamma^\infty(\Anti^k TE^*)$ on $E^*$ we consider $X_t =
    \E^{kt} \Phi_t^* X \in \Gamma^\infty(\Anti^k TE^*)$. A simple
    argument (say in local coordinates) shows that the limit $t
    \longrightarrow - \infty$ of $X_t$ exists and is a vertical
    lift. We def\/ine
    \[
    \Unlift(X) = \lim_{t \rightarrow -\infty} X_t.
    \]
    For a vertical lift $\omega^\ver$ we clearly have
    $\Unlift(\omega^\ver) = \omega^\ver$ and $\Unlift$ is a projection
    onto the vertical lifts $\Gamma^\infty(\Anti^kE^*)^\ver \subseteq
    \Gamma^\infty(\Anti^k TE^*)$. Since $\theta_E$ is homogeneous of
    degree $-1$, i.e.\ $\E^{-t}\Phi_t^* \theta_E = \theta_E$, and since
    the Schouten bracket is natural with respect to pull-backs we
    f\/inally obtain by continuity
    \[
    \Unlift(\D_{\theta_E} X)
    = \lim_{t \rightarrow -\infty} \E^{k-1}
    \Phi_t^* \Schouten{\theta_E, X}
    = \Schouten{\theta_E, \lim_{t \rightarrow -\infty} \E^{k} \Phi_t^*X}
    = \D_{\theta_E} \Unlift(X).
    \]
    Thus the projection $\Unlift$ descends to cohomology from which
    the injectivity of \eqref{eq:VerticalLiftCohomology} follows
    immediately.
\end{proof}
\begin{Remark}
    \label{remark:Exercise}
    In particular, this theorem shows that the scalar
    Chevalley--Eilenberg cohomology of a Lie algebra $\mathfrak{g}$ is
    injected into the Poisson cohomology of the linear Poisson
    structure on the dual $\mathfrak{g}^*$ of a Lie algebra
    $\mathfrak{g}$, see
    \cite[Exercise~72]{cannasdasilva.weinstein:1999a}. Note however,
    that e.g.\ for Abelian Lie algebras this map is far from being
    surjective.
\end{Remark}

Recall that the modular vector f\/ield $\Delta_\Omega \in
\Gamma^\infty(T E^*)$ with respect to
some positive density $\Omega \in \Gamma^\infty(|\Anti^{n+N}|T^*E^*)$
is def\/ined by
\begin{equation*}
    \Delta_\Omega(f) \Omega = \Lie_{X_f} \Omega
\end{equation*}
and gives a Poisson vector f\/ield $\Delta_\Omega$ which depends on
$\Omega$ only via Hamiltonian vector f\/ields. Thus the \emph{modular
  class} $[\Delta_\Omega] \in \mathrm{H}^1_{\theta_E}(E^*)$ is
well def\/ined, see e.g.\ \cite{evens.lu.weinstein:1999a,
  weinstein:1997a}. Since $\Pol^0(E^*)$ and $\Pol^1(E^*)$ generate
$\Pol^\bullet(E^*)$ we can conclude from \eqref{eq:LieXJsmunu} and
\eqref{eq:LieXumunu}
\begin{equation*}
    \Delta_{\mu \otimes \nu} = \left(\tr\ad\right)^\ver,
\end{equation*}
whence $[\Delta_{\mu \otimes \nu}]$ is precisely the image of the
modular class under~\eqref{eq:VerticalLiftCohomology}.  From the
injectivity of~\eqref{eq:VerticalLiftCohomology} we immediately obtain
the following corollary, see also \cite[Section~7]{weinstein:1997a}:
\begin{Corollary}
    \label{corollary:Unimodular}
    Let $E$ be a Lie algebroid. Then the following statements are
    equivalent:
    \begin{enumerate}
    \itemsep=0pt
    \item\label{item:Eunimodular} $E$ is unimodular in the sense of
        Lie algebroids.
    \item\label{item:Estarunimodular} $(E^*,\theta_E)$ is unimodular
        in the sense of Poisson manifolds.
    \item\label{item:ConstantOmega} There exists a positive density
        $\Omega \in \Gamma^\infty(|\Anti^{n+N}|T^*E^*)$ on $E^*$ which
        is constant along the fibres and satisfies $\Lie_{X_f} \Omega
        = 0$ for all $f \in C^\infty(E^*)$.
    \end{enumerate}
\end{Corollary}

Now we assume that $E$ and hence $E^*$ are unimodular and we choose an
appropriate positive density $\Omega \in
\Gamma^\infty(|\Anti^{n+N}|T^*E^*)$ on $E^*$ which is constant along
the f\/ibres with $\Delta_\Omega = 0$. Then we can proceed literally as
in \cite[Section~8]{bordemann.neumaier.waldmann:1999a} and
\cite{bieliavsky.bordemann.gutt.waldmann:2003a} to show that the
functional
\begin{equation}
    \label{eq:PoissonTrace}
    \tr(f) = \int_{E^*} f \Omega
\end{equation}
for $f \in C^\infty_0(M)[[\lambda]]$ def\/ines a
$\mathbb{C}[[\lambda]]$-linear trace for any \emph{homogeneous} star
product on $E^*$. Thus let $\star = \sum_{r=0}^\infty \lambda^r C_r$
be a homogeneous star product quantizing the Poisson bracket
$\{\cdot,\cdot\}_E$ of~$E^*$.  The f\/irst lemma uses the fact that
$\Omega$ is constant along the f\/ibre directions and is shown by
partial integration.
\begin{Lemma}
    \label{lemma:HomogeneousDiffOp}
    Let $D: C^\infty(E^*) \longrightarrow C^\infty(E^*)$ be a
    homogeneous differential operator of homogeneity $-r$ with $r \ge
    1$, i.e.\ $[\Lie_\xi, D] = -r D$. Then $\tr(D(f)) = 0$ for all $f
    \in C^\infty_0(E^*)$.
\end{Lemma}
\begin{Lemma}
    \label{lemma:TraceStarHomogeneous}
    Let $f \in \Pol^k(E^*)$ and $g \in C^\infty_0(E^*)$. Then
    \begin{equation*}
        \tr(f \star g)
        = \tr\left(\sum_{r=0}^k \lambda^r C_r(f, g)\right)
    \end{equation*}
    and analogously for $\tr(g \star f)$.
\end{Lemma}
\begin{proof}
    This follows from the preceding lemma as $C_r(f, \cdot)$ and
    $C_r(\cdot, f)$ are homogeneous dif\/fe\-ren\-tial operators of
    homogeneity $k-r$.
\end{proof}
\begin{Lemma}
    \label{lemma:TraceIsTraceOnPol}
    Let $f \in \Pol^\bullet(E^*)$ and $g \in C^\infty_0(E^*)$. Then
    \begin{equation}
        \label{eq:TraceIsTraceOnPol}
        \tr(f \star g - g \star f) = 0.
    \end{equation}
\end{Lemma}
\begin{proof}
    Let $f$ be homogeneous of degree $k$. For $k=0$ the statement is
    fulf\/illed by Lemma~\ref{lemma:TraceStarHomogeneous}. For $k=1$ it
    is fulf\/illed by Lemma~\ref{lemma:TraceStarHomogeneous} and the
    unimodularity condition. For $k > 1$ it is suf\/f\/icient to consider
    $f = X_1 \star \cdots \star X_k$ with $X_i$ homogeneous of degree
    $0$ or $1$ according to
    Proposition~\ref{proposition:HomogeneousStarProducts}~\refitem{item:PolynomialGenerated}.
    From
    \[
    [X_1 \star \cdots \star X_k, g]_\star =
    [X_1, X_2 \star \cdots \star X_k \star g]_\star
    +
    [X_2 \star \cdots \star X_k, g \star X_1]_\star
    \]
    we conclude \eqref{eq:TraceIsTraceOnPol} by induction on $k$.
\end{proof}
\begin{Theorem}
    \label{theorem:UnimodularTrace}
    Let $\star$ be a homogeneous star product on $E^*$ and let $\tr$
    be the Poisson trace as in~\eqref{eq:PoissonTrace}. Then $\tr$ is
    a trace functional with respect to $\star$, too.
\end{Theorem}
\begin{proof}
    By the usual Stone--Weierstra\ss{} argument and the continuity of the
    integration the theorem follows from the last lemma.
\end{proof}
\begin{Remark}
    Of course there are much more general statements on the traces of
    star products available like the cyclic formality theorem, see the
    discussion in \cite{felder.shoikhet:2000a}. In the case of complex
    Lie algebroids, one may also consult \cite{chemla:1999a}. Note
    however, that in the particular situation we are discussing, the
    proof is elementary.
\end{Remark}

\subsection*{Acknowledgements}

We would like to thank Janusz Grabowski,
Simone Gutt, Yvette Kosmann-Schwarzbach and Alan Weinstein for
valuable discussions and remarks. Moreover, we thank the referees for
many interesting suggestions and remarks.

\section*{Editorial comments}


For the benef\/it of the readers we reproduce here
verbatim the sug\-gestion of the referee mentioned in Remark~\ref{remark:Extension}:

``We believe that most of the results presented in the manuscript under consideration can be
recovered in a very simpler way. Namely, to any ${\mathbb C}[[\hbar]]$-valued 2-cocycle $B_\hbar$ as above one can
associate a central extension $\widetilde{E}_\hbar$ of $E[[\hbar]]$ (def\/ined over ${\mathbb C}[[\hbar]]$) and can consider its universal
enveloping algebra $U(\widetilde{E}_\hbar)$. But we instead consider the $\hbar$-universal enveloping algebra
$U_\hbar(\widetilde{E}_\hbar)$ (i.e.\ we put an $\hbar$ in front of commutators in the def\/ining relations).

Finally, take the quotient $A_\hbar$ of $U_\hbar(\widetilde{E}_\hbar)$ by $c = \hbar$, where $c$ is the generator of the
one-dimensional extension. One can prove that $A_\hbar$ is a topologically free ${\mathbb C}[[\hbar]]$-module
(i.e.\ $A_\hbar \cong \Gamma(M, S(E))[[\hbar]]$), and is a quantization of the $B$-twisted Poisson structure on
functions polynomial in the f\/ibers on $E^*$. Two such quantized algebras are obviously isomorphic
if they come from cohomologous cocycles, and the homogeneity property follows
easily from this simple construction in the case when $B_\hbar = B$.

The proof that this quantization extends to the full algebra of smooth functions on $E^*$
follows from the already known case when $B = 0$.

\medskip

\noindent
{\bf Remark.} In the case of a Lie algebra this construction is very well-known (it has been
used by Drinfel'd to quantize triangular $r$-matrices, and it is quite standard in the theory of
quantum integrable systems). In the case of the cotangent bundle these are called twisted
dif\/ferential operators.''

\pdfbookmark[1]{References}{ref}
\LastPageEnding

\end{document}